\documentclass[12pt]{amsart}

\usepackage[T1]{fontenc}
\usepackage{amsfonts,amsmath,amssymb,amsthm,tikz}
\usepackage[shortlabels]{enumitem}
 \usepackage{lineno}
\usetikzlibrary{decorations.pathreplacing}
\usepackage{float} 

\usepackage[pagebackref]{hyperref} 

\usepackage{multicol}
\usepackage{color}

\usepackage{epsfig}
 \usepackage[abbrev, nobysame]{amsrefs}
 
\usepackage[latin1]{inputenc} 

\usepackage{lineno} 

\makeatletter
\@namedef{subjclassname@2020}{%
  \textup{2020} Mathematics Subject Classification}
\makeatother

\newtheorem{Lem}{Lemma}[section]
\newtheorem{The}[Lem]{Theorem}

\newtheorem{Cor}[Lem]{Corollary}
\newtheorem{Pro}[Lem]{Proposition}

\theoremstyle{definition}

\newtheorem{Defi}[Lem]{Definition}
\newtheorem{Exam}[Lem]{Example}
\newtheorem{Not}[Lem]{Notation}
\newtheorem{Rem}[Lem]{Remark}

\newtheorem*{notation}{Notation}

\def\a{\alpha}

\def\b{\beta}
\def\C{\mathbb C}
\def\D{\Delta}
\def\d{\delta}

\def\f{\phi}

\def\g{\gamma}

\def\ll{\lambda}
\def\N{\mathbb N}

\def\ord{\mbox{\rm ord}}

\def\O{\mbox{\rm orb}}
\def\p{\pi}

\def\Q{\mathbb Q}
\def\R{\mathbb R}
\def\r{\rho}

\def\s{\sigma}
\def\t{\tau}

\def\y{\wedge}
\def\Z{\mathbb Z}
\def\z{\zeta}

\begin{document}

\title{Logarithmic Jacobian ideals, quasi-ordinary hypersurfaces and equisingularity}
\author{Pedro Daniel Gonz\'alez P\'erez}

\address{Instituto de Matem\'atica Interdisciplinar (IMI) \\ Departamento de \'Algebra, Geometr\'{\i}a y Topolog\'{\i}a \\
Facultad de Ciencias Matem\'aticas \\ Universidad Complutense de Madrid\\
Plaza de las Ciencias 3. 28040. \\
Madrid. Spain}

\email{pgonzalez@mat.ucm.es}

\keywords{Jacobian ideal, logarithmic Jacobian ideal, toric
singularities, quasi-ordinary singularities, Nash modification, equisingularity}

\date{\today}
\subjclass[2000]{Primary 14M25; Secondary 32S25}

  \begin{abstract}

  Let $(S, 0) \subset (\C^{d+1},0)$ be an irreducible germ of hypersurface. 
The germ $(S,0)$  is quasi-ordinary if $(S,0)$ has a finite projection to 
   $(\C^d,0)$ which is unramified
outside the coordinate hyperplanes.  This implies that the normalization of $S$ is a toric singularity.
One has also a monomial variety associated to $S$, which is a toric singularity with the same normalization, and  with
possibly higher embedding dimension. 
Since $(S,0)$ is quasi-ordinary, the extension of the Jacobian ideal of $S$ to the local ring of its normalization is a monomial ideal. 
We describe this monomial ideal by comparing it with the {\em logarithmic Jacobian ideals} of $S$ and of its associated monomial variety and we give some applications.
\end{abstract}

\maketitle

{\em \hfill Dedicated to the memory of Arkadiusz P\l oski. \smallskip}


\section*{Introduction}

An equidimensional germ $(S, 0)$ of complex analytic variety of
dimension $d$ is quasi-ordinary if there exists a {\em quasi-ordinary projection}, that is, a finite map $\p:
(S, 0) \rightarrow (\C,^d  0)$, which is unramified outside a
normal crossing divisor in  $(\C^d , 0)$. This class of
singularities plays an important role in the classical approach to
study singularities by the Jung's method (see
\cite{Lipman2} and \cite{Abhyankar}).

In this paper we suppose that $(S, 0)$ is an analytically irreducible hypersurface germ
of dimension $d \geq 1$. 
The quasi-ordinary hypersurface $(S,0)$ may be defined by a quasi-ordinary polynomial 
$f \in \C \{ x_1, \dots, x_d \} [y]$, that is, a Weierstrass polynomial  such that its discriminant 
$\D_y f$ is of the form a monomial times a unit in the ring $\C \{ x_1, \dots, x_d \}$.
The germ $S$ has
certain fractional power series parametrizations $y = \z(x_1^{1/n}, \dots, x_n^{1/n})$ for $n = \deg_y f$,
which generalize the
Newton-Puiseux series of plane curve singularities. 
These fractional power series have a finite set of 
 {\em characteristic or distinguished monomials}, 
 which generalize the
characteristic exponents in the case of plane branches. 
Joseph Lipman and Yih-Nan Gau proved that these monomials classify the embedded
topological type of $(S, 0) \subset (\C^{d+1}, 0)$     (see
\cite{Lipman2, Gau}).

As in the case of plane branches, one has a notion of semigroup $\Gamma$ 
associated to the  quasi-ordinary singularity $(S, 0)$.
The semigroup $\Gamma$ determines and it is determined by the characteristic monomials, and it  is independent of the choice of
parametrization of $(S, 0)$, see \cite{PPP04, GP1, KM}.
The semigroup $\Gamma$ defines an
 affine toric variety $Z^\Gamma$, which is a complete intersection of dimension $d$, 
called the \textit{monomial variety} associated with $S$.


Quasi-ordinary hypersurface singularities are well adapted to toric geometry methods: for instance the normalization
$(Z, 0) \rightarrow (S, 0)$ is a toric singularity determined
by the characteristic monomials (see \cite{GP4,
PPP04}). 
The composition $(Z,0) \to (\C^d, 0)$ of the normalization with the quasi-ordinary projection 
is also a quasi-ordinary projection.
Recall that the singular locus of $S$ is defined by the Jacobian ideal $\mathrm{Jac}(S)$ of $S$. 
It follows that the pullback of the singular locus of $S$ by the normalization map is contained in the complement of the torus of $Z$.
This implies that the ideal $\mathrm{Jac}(S) \mathcal{O}_Z$ is a monomial ideal with respect to the toric structure
of the ring $\mathcal{O}_Z$. The main goal of this paper is to provide a description of the monomial
ideal $\mathrm{Jac}(S) \mathcal{O}_Z$. 

This goal is related to the study of the normalized Nash modification. 
Recall that the Nash modification of an algebraic (or complex analytic)
variety is a canonical morphism which replaces each singular point
by the limiting positions of tangent spaces at non-singular
points. 
Since $S $ is a complete intersection, the Nash
modification of $S$ is isomorphic to the blow up of the Jacobian
ideal in $\mathcal {O}_{S}$ (see \cite{Nobile}
and also \cite{Piene} for a similar result in the case of Gorenstein singularities). Then, the 
normalized Nash modification of $S$, which is the Nash modification followed by the normalization, 
 is isomorphic  to the composition
of the normalization map of $S$ with the normalized blow up of
${\textrm{Jac}(S)} \mathcal{O}_{Z} $. This follows from the
properties of normalized blow ups of ideals 
(see \cite{LT}, Propositions 3.2 and 3.3). 

The Nash
modification of an affine toric variety defined over an algebraically closed field
of characteristic zero is isomorphic to the blow up of a monomial ideal, called the {\em
logarithmic Jacobian ideal} (see 
 \cite{GN, LR, T3}). 
 We define the  logarithmic Jacobian ideal $\mathcal{J}_{\mathrm{log}}(V)$ of a
$d$-dimensional variety $V$ with toric normalization $\overline{V}$ as
the image of the module of differentials $\Omega_{V}^d$ by the
composite of the maps:
\[
\Omega_{V}^d \longrightarrow  \Omega_{\overline{V}}^d \longrightarrow
\Omega_{\overline{V}}^d (\log D) \longrightarrow \mathcal{O}_{\overline{V}},
\]
where $\Omega_{\overline{V}}^d (\log D)$ is the module  of $d$-forms
with logarithmic poles on the complement $D$ of the torus of
$\overline{V}$, the first map is induced by the normalization, the
second map is the canonical one from the toric structure of $\overline{V}$
and the third map is an isomorphism.


We describe the Jacobian ideal of an affine toric variety  in Proposition \ref{Jacobian-ideal-toric}.
By Remark  \ref{rem:pzero} this result implies, in the case $V$ is the affine toric variety $Z^\Gamma$ associated to the semigroup $\Gamma$ of a quasi-ordinary hypersurface singularity $(S, 0)$, that
\begin{equation} \label{eq:inc2}
X^{\g_0} \,   \mathcal{J}_{\mathrm{log}}(Z^\Gamma)  =   {\textrm{Jac} (Z^\Gamma)},
\end{equation}
where  $\g_0$ is the minimal Frobenius vector of the semigroup $\Gamma$ (see Proposition \ref{def: Frobenius}).


It is a natural to ask if the relation \eqref{eq:inc2} can be extended to the quasi-ordinary hypersurface germ $(S,0)$.
The ideal  $\mathcal{J}_{\mathrm{log}}(S)$ and the Frobenius vector $\g_0$ are determined by 
the characteristic monomials, hence by the embedded topological type of $(S,0) \subset (\C^{d+1}, 0)$. The ideal $\mathcal{J}_{\mathrm{log}}(S)$ was described by  Helena Cobo and the author in  \cite{pseries-qo}. 
Our main results are Theorem \ref{th:inclusion}, which shows  the inclusion
\begin{equation} \label{eq:inc}
X^{\g_0} \,   \mathcal{J}_{\mathrm{log}}(S) \mathcal{O}_{Z} \, \subseteq   \, {\textrm{Jac} (S)} \mathcal{O}_{Z},
\end{equation}
and  Corollary \ref{cor:inclusion1}, which states that equality holds in \eqref{eq:inc}
in the two dimensional case.


We will apply a particular
re-embedding of the germ $S$ in an affine space of dimension $d+ g$, where $g$ is the number of characteristic monomials.  
Using this embedding we can build a  deformation with fibers $S_p$,  
for $p \in Z^\Gamma$, such that if $p$ belongs to the torus  then
$S_p$ is isomorphic to the germ $S$, 
while if $p$ is equal to the origin $0$  then  the special fiber $S_0$ is isomorphic to $Z^\Gamma$.
In both cases the fibers have the same normalization
$\overline{S}_p= Z$. There is a toric modification of $\C^{d+g}$ which provides an embedded 
resolution of $S_p$, for $p$ in the torus or $p=0$ (see  \cite{GP4}).
In the case of plane branches this result
was proved by Rebeca Goldin and Bernard
Teissier in  \cite{Goldin}. 


In order to describe the Jacobian ideal $\mathrm{Jac}(S_p)$, 
we introduce a  matrix $R(p)$ in terms of 
the Jacobian matrix of the functions defining the fibers of this deformation. 
The matrix $R(p)$ specializes when $p =0$  to a matrix with integer coefficients encoding the 
relations among the
generators of the semigroup $\Gamma$.


As an application,  we give a method to compute the number 
 $\overline{\nu}_{\textrm{Jac} (S)} ( \phi) $  associated to  a function $\phi \in \mathcal{O}_S$, when the dimension of $S$ is equal to two. The number $\overline{\nu}_{\textrm{Jac} (S)} ( \phi) $ is 
 studied by Monique Lejeune-Jalabert and Bernard Teissier in a more general setting in \cite{LT}.
 The number
 $\overline{\nu}_{\textrm{Jac} (S)} ( \phi) $ can be expressed in terms of $\phi$ and
the divisorial valuations associated with the
irreducible components of the pull-back of the Jacobian ideal of
$S$  in the normalized Nash modification (see \cite{LT}).
By \cite[Section 6]{LT}, the number $({\overline{\nu}_{\textrm{Jac} (S)} ( \phi)})^{-1}$ 
may be seen as a {\em \L ojasiewicz exponent} associated to the ideal ${\textrm{Jac} (S)}$ 
and the function $\phi$ in a compact neighborhood of the $0 \in S$. 
The study of \L ojasiewicz exponents has been one of the research interests of  Arkadiusz P\l oski, and he has published many articles on this topic, either alone or with his collaborators 
 (see for instance 
 \cite{MR817547,
 MR1101851,
 MR786190,
 MR805025,
 MR1350859,
 MR2553365,
 MR2882792,
  MR2148949,
   MR2515408,
    MR3538652}).


The results and proofs of this paper hold in the category of
complex analytic spaces, and also in the algebroid category over an
algebraically closed field  of zero characteristic.

%
The structure of the paper is as follows.  In Section
\ref{nash-tor2} we introduce the basic notions of toric geometry, we review the construction of the Nash
modification of an affine toric variety and we give a combinatorial
description of its Jacobian ideal.
In Section \ref{qo} we review the basic properties of quasi-ordinary hypersurface singularities 
and their associated semigroups. 
In Section \ref{rel-Jacobian} we recall the construction of the logarithmic Jacobian ideals of a toric and of a quasi-ordinary hypersurface. In Section \ref{Sec:JI} we describe the deformation of $S$ which specializes to the monomial variety $Z^\Gamma$ and then we prove the main results about the relations between the Jacobian ideal and the logarithmic Jacobian ideals. Finally, in Section \ref{sec:app} we give an application to the computation of the number
 $\overline{\nu}_{\mathrm{Jac}(S)} (\phi)$ associated to a function $\phi \in \mathcal{O}_S$.

\begin{notation} 
$\,$ 

\noindent 
$\bullet$ The result of    erasing the term  $r_j$ in the sequence   $r_1, \dots, r_s$ is denoted by
$r_1, \dots ,\hat{r}_j, \dots, r_s$.

\noindent 
$\bullet$ Let $A= (a_{j}^{i})^{i=1, \dots, k} _{j=1, \dots, l}$ be a matrix
with $k$ rows and $l$ columns, and with coefficients in a ring. If $1\leq i_1 < \cdots < i_s \leq k$
(resp. if  $1 \leq j_1 < \dots < j_r \leq l$)  we denote by
$A^{i_1, \dots, i_s}$ (resp. by $A_{j_1, \dots, j_r}$) the
submatrix
\[ A^{i_1, \dots, i_s} := (a^{i}_{j})_{j\in \{1, \dots
l\}} ^{i \in \{1, \dots, \hat{i}_1, \dots, \hat{i}_s,\dots,  k \} }
\quad (\textrm{resp. } A_{j_1, \dots, j_r} := (a^{i}_{j})^{i \in
\{1, \dots, k \}}_{ j\in \{ 1, \dots, \hat{j}_1, \dots, \hat{j_r}, \dots,  l \} } ).
\] 
\noindent 
$\bullet$ If $A$ is a square matrix we denote by $|A|$ its
determinant. If $A$ is a matrix we denote its rank by $\mathrm{rk} A$.

\noindent 
$\bullet$ If  $\g, \g' \in \R^d$ we set 
\begin{equation} \label{eq-order-def}
\g \prec \g'  \quad \mbox{ iff }  \quad \g' - \g \in \R^d_{\geq 0}  \quad 
\mbox{ and  } \quad  \g \ne \g'.
\end{equation}
We write $\g \preceq \g'$ iff $\g' - \g \in  \R^d_{\geq 0}$.

\end{notation}

\section{On the Nash modification of an affine toric variety}
\label{nash-tor2}

\subsection{Equivariant embedding of an affine toric variety}     \label{sec-tor}
We give first some basic definitions
and notations (see \cite{Fulton,St, Cox}).


 If $N \cong \Z^{d}$ is a lattice we denote by $M$ its
dual lattice, by $N_\R$ (resp. $N_\Q$) the vector space spanned by
$N$ over the field $\R$ (resp. $\Q$). We use similar notations for a
lattice homomorphism $\f$ and the associated vector space
homomorphism $\f_\Q$ (or $\f_\R$).
We denote by \[ \langle \,  \,  , \,  \, \rangle: N\times M \to \Z,  \quad (u,v) \mapsto \langle u, v \rangle,\] the dual pairing between 
the lattices $N$ and $M$.
In what follows a {\it cone} means a
{\it rational convex polyhedral cone}: the set of   nonnegative
linear combinations of vectors $v_1, \dots, v_r \in N$. The cone
$\sigma$ is {\it strictly convex} if it contains no lines.
The {\it dual} cone  $\sigma^\vee$
of $\sigma$ is the set $ \{ w  \in M_\R \mid \langle w, u \rangle
\geq 0 \}$.
We denote by $\mathrm{int}({\s})$ the {\em relative interior} of the cone $\s$.


Let  $ \Lambda $ be a sub-semigroup of finite type of the lattice $M$.
We will assume that $\Lambda$ generates $M$ as a group, that is,  $M = \Z \Lambda$. The
$\C$-algebra of the semigroup $\C [ \Lambda ] := \{ \sum_{\mathrm{finite}} a_\ll
{X}^\ll \mid a_\ll \in \C \}$ is of finite type and corresponds to
the affine toric
  variety
of the form $Z^{\Lambda} = \makebox{Spec} \, \C [ \Lambda ]$. The
torus $Z^{M}$ is an open dense subset of $Z^\Lambda$, which acts
on  $Z^\Lambda$ in a way which extends the group multiplication on the torus.
A closed point $p \in Z^\Lambda$ determines and it is determined by the semigroup homomorphism 
$g_p : \Lambda \to (\C, \cdot)$, such that $g_p (m) = X^m(p)$ for $m \in \Lambda$, where $(\C, \cdot)$ represents the semigroup defined by $\C$ with the multiplication. For instance, the unit point $\underline{1}$ of the torus $Z^{M}$ corresponds to the constant homomorphism 
$g_{\underline{1}}$ defined by 
$g_{\underline{1}} (m) = 1$, for all $m \in M$. 
We have a bijection $\t \mapsto \O _\t$ between the faces
of $\s$ and the orbits of the torus action which reverses the
inclusion of the closures.     Any affine toric variety is defined
by a sub-semigroup of this form.   


The
cone $\s^\vee := \R_{\geq 0} \Lambda \subset M_\R$ is rational for
the lattice $M$, hence the semigroup $ \s^\vee \cap M$ is of
finite type. The map of $\C$-algebras $\C[\Lambda] \rightarrow \C[
\s^\vee \cap M]$ induced by
 $\Lambda \hookrightarrow \s^\vee \cap
M $,  corresponds geometrically to the normalization map
$Z^{\s^\vee \cap M} \rightarrow Z^{{\Lambda}}$.


From now on we assume that the cone $\s^\vee := \R_{\geq 0}
\Lambda$ is strictly convex of dimension $d$. In this situation
the zero dimensional orbit of the torus action is reduced to the origin of $Z^\Lambda$, which is 
the point $ 0 \in Z^{\Lambda}$ defined by the maximal ideal $
\mathfrak{m}_\Lambda := (\Lambda - \{ 0 \} ) \C [ \Lambda ]$. The origin $0$ of $Z^\Lambda$ corresponds to 
the homomorphism of semigroups $g_0: \Lambda \to (\C^*,0)$ such that 
$g_0 (m) =0$ for all $m \in \Lambda \setminus \{0\}$ and $g_0 (0) = 1$. 


We denote by  $ \C[[\Lambda] ]$ the ring
of formal power series with coefficients in $\C$ and exponents in
the semigroup $\Lambda$.   
We denote by $\C \{ \Lambda \}$ the local ring of germs of
holomorphic functions of $Z^\Lambda$ at the origin. 
Its completion with respect to its maximal ideal is 
the ring  $\C [[ \Lambda ]]  = \{ \sum a_\ll
{X}^\ll \mid a_\ll \in \C \}$ 
of formal power series with exponents in the semigroup $\Lambda$ (see \cite{GP4}).


 Let $\a_{1}, \dots,
\a_{m} \in \Lambda \setminus \{ 0 \}$ be a set of generators of $\Lambda$. 
Denote by $e_1, \dots, e_m$ the canonical basis of $\Z^m$.
We have an exact sequence of lattice homomorphisms 
\begin{equation} \label{eq:sequence}
0 \longrightarrow \mathcal{L}  \longrightarrow \Z^m \stackrel{\psi}{\longrightarrow} M  \longrightarrow 0 ,
\end{equation}
where $\psi( e_j) = \a_j$, for $j = 1, \dots , m$, and $\mathcal{L}$ is the kernel of $\psi$. 
The surjective 
homomorphism of semigroups
\begin{equation} \label{eq:presentation}
\Z^m_{\geq 0} \rightarrow \Lambda,  \quad e_j \mapsto
\a_{j}, \quad j =1, \dots, m,
\end{equation}
corresponds to an equivariant embedding $Z^\Lambda \hookrightarrow
\C^m$ given by
\begin{equation} \label{iota}
 \psi_*: \C[U_1, \dots, U_m]
\rightarrow \C[\Lambda], \quad U_i \mapsto {X}^{\a_{j}} , \quad
\mbox{ for } j=1, \dots, m. \end{equation} 

The {\em support} of a vector $n = \sum_{j=1}^m n_j e_j\in \Z^m$ is the set 
$\mathrm{supp} (n) = \{ j \in \{1, \dots, n \}
\mid n_j \ne 0 \}$. 
Let $n \in \Z^m$ be a nonzero vector in
$ \mbox{\rm ker} (\psi)$. 
The vector $n$ decomposes in a unique way in 
the form $n = n^+ -n^-$, where $n^+, n^- \in \Z^m_{\geq 0}$ and 
have disjoint support. 
 Notice that the vectors 
$n^+$ and $n^- $ are nonzero,
otherwise we would get a non-trivial linear combination $0 = \sum_j n_j \a_j$ in $M$,
which contradicts the strict convexity of the cone $\s^\vee$.
 Since $n \in \ker \psi$, we have the relation
\begin{equation} \label{relation}
\psi(n^+) = \sum_{i=1}^m n_i^+ \a_{i}=  \sum_{i=1}^m n_i^- \a_{i} = \psi(n^-).
\end{equation}
There exists 
a sequence of vectors
$n_{1}, \dots, n_{p} \in \Z^m
$ such 
that the ideal $\mathcal{I}^\Lambda$ of the embedding $Z^\Lambda \hookrightarrow \C^m$
is generated by the binomials $h_j := U^{n_{j}^+} - U^{n_{j}^-}$, for $j=1, \dots,
p$.  
The vectors $n_{1}, \dots, n_{p}$ generate $\mbox{\rm ker} (\psi)$.
The semigroup $\Lambda$ is a \textit{complete intersection} 
if the minimal number of binomial generators
of  the ideal $\mathcal{I}^\Lambda$ is equal to $m- d$.
This means that $Z^\Lambda$ is a \textit{complete intersection}.


We have an exact sequence
\begin{equation} \label{sequence}
 \Z^p
\stackrel{\f}{\longrightarrow} \Z^{m}
\stackrel{\psi}{\longrightarrow} M \longrightarrow 0,
\end{equation}
in which the first map ${\f}$ applies the
$j$th canonical vector to the 
vector $n_j$, for $j =1, \dots, p$. 
Write $n_{i}= \sum_{j=1}^m \ell_j^{(i)} e_i $, for $i =1, \dots, p$.
The 
{\em matrix
 of
relations} 
associated with $n_{1}, \dots, n_{p}$ is 
\begin{equation} \label{rel2}
R:= \left(
\begin{array}{ccc}
\ell^{(1)}_{1} &  \cdots & \ell^{(1)}_{m}
\\
\cdots &  \cdots & \cdots
\\
\ell^{(p)}_{1} &  \cdots & \ell^{(p)}_{m}
\end{array}
\right).
\end{equation}
The matrix $R$ has rank $m - d \leq p$.
Notice that we can always relabel the vectors $n_1, \dots, n_p$ in such a way
that the submatrix $R^{m-d+1, \dots, p}$ 
of $R$
 is of rank $m-d$.

\begin{Lem}  \label{linear} (see \cite[Proposition 60]{T3}, \cite[Lemma 1.8]{MR4792759}) 
If the matrix $R^{m-d+1, \dots, p}$ 
 is of rank $m-d$ then for any sequence
$1\leq j_1 < \dots < j_d \leq m$ we have 
\[
\a_{j_1} \y \cdots  \y \a_{j_d} = 0 \mbox{ if and only if } |R^{m-d+1,
\dots, p} _{j_1,\dots, j_d}| = 0.\]
\end{Lem}
\begin{proof}
We consider the exact sequence induced by
(\ref{sequence}) over the field  $\Q$. The linear subspace $L:=
\mbox{\rm Im} (\f_\Q)$, is generated by the images of the first
$m-d$ canonical vectors by hypothesis. The matrix associated to
the  restriction $\f_\Q'$ of the linear map   $\f_\Q$ to the subspace
spanned by the first $m-d$ canonical vectors, with respect to the
canonical basis,  is equal to the transpose of  $R^{m-d+1, \dots,
d}$.

The matrix $P$ associated to $\psi_\Q$ with respect to the canonical
basis of $\Q^{m}$ and a fixed basis of $M_\Q$ has $i$th column
the coordinates of $\a_i \in M_\Q $, for $i=1,\dots, m$. By
dualizing the sequence (\ref{sequence}) we notice that the
$j$th row of the matrix $P$ defines the coordinates of a linear form $w_j \in (\Q^{m})^*$, with respect to the dual basis, in
such a way that the vector subspace $L$ is also obtained by the
intersection of the kernels of $w_j$, for $j =1,\dots, d$.

We have exhibited the $m-d$-dimensional subspace $L$ 
 by equations and by generators. We deduce the
assertion as a consequence of the classical relations of the
Grassmann coordinates of a linear subspace and its dual Grassmann
coordinates applied to $L$ (see \cite{HP}, VII, \S 3, Theorem
1).
\end{proof}


A vector $\a_0 \in M$ is a \textit{Frobenius vector} of the semigroup $\Lambda$ if 
for every  vector $\a \in \mathrm{int} ( \s^\vee) \cap M$,
the sum $\a_0 +\a$ 
belongs to $\Lambda$ (see \cite{AGO15}). This notion generalizes the Frobenius number of a numerical semigroup, in the case
$M= \Z$. For instance, if $e_1, e_2$ is the canonical basis of $\Z^2$ and then $-e_1 -e_2$ is a Frobenious vector of the semigroup $\Z^2_{\geq 0}$.


We say that $\Lambda$ is a \textit{complete intersection semigroup} if the ideal $\mathcal{I}^\Lambda$ is minimally generated by $m - d$  binomials.  If $\Lambda$ is a complete intersection semigroup then $\Lambda$ has a \textit{minimal} Frobenius vector $\a_0$ (see \cite[Th. 3.3]{AGO15}). The minimality means
that 
if $\a_0'$ is a Frobenius vector then $\a_0 ' \in \a_0 + \s^\vee$.

\subsection{On the Nash modification of a toric singularity}       \label{nash-tor}

In this section we review the description of the Nash modification of a complex affine toric variety $Z^\Lambda$ following \cite{T3}. As an application, we obtain a description of the generators of the Jacobian ideal $\mathrm{Jac} (Z^\Lambda)$
in Proposition \ref{Jacobian-ideal-toric}.

\begin{Not} 
We consider a semigroup $\Lambda$ as in Section
\ref{sec-tor}.
We denote by $J$ (resp. $\tilde{J}$) the matrix with coefficients
in $\C[\Lambda]$ defined from the Jacobian matrix of $(h_1, \dots,
h_p)$ by
\[
J := \left( \psi_* (\frac{{\partial} h_i}{\partial U_j})
\right)_{j=1, \dots, m}^{i=1, \dots, p}, \quad \left(\textrm{resp. } 
\tilde{J} := \left( \psi_* (U_j  \cdot \frac{{\partial} h_i}{\partial
U_j}) \right)_{j=1, \dots, m}^{i=1, \dots, p} \right),
\]
where we recall that $\psi_*$ is defined by \eqref{iota}
\end{Not}

\begin{Pro}      \label{nash-toric}
(see
  \cite{GN, T3, LR})
The Nash modification of $Z^{\Lambda} $ is the blow up of the
ideal  $ \mathcal{J}_\mathrm{log} (Z^\Lambda)$ of $\C[\Lambda]$ generated by the images of the products
$U_{j_1} \dots U_{j_d}$ such that $\a_{j_1} \y \dots \y \a_{j_d}
\ne 0$, i.e., by the monomial ideal
\begin{equation}   \label{jota-lambda}
 \mathcal{J}_\mathrm{log} (Z^\Lambda) := (\{ X^{\a_{j_1}  +
  \cdots+  \a_{j_d}}  \, \mid \, 1 \leq j_1 < \cdots < j_d
  \leq m , \, \a_{j_1}  \y
  \cdots \y \a_{j_d} \ne 0 \}).
\end{equation}
\end{Pro}
\begin{proof}
 We suppose first that the toric singularity
$Z^\Lambda$ is a complete intersection. This means that  we can assume $p = m-d$,
and the map $\f$ in the sequence (\ref{sequence}) is injective.


Nobile proved that the Nash modification of a complete
intersection $X$ is the blow up of the image of the Jacobian ideal
of $X$ in its coordinate ring ${\mathcal O}_X$  (see \cite{Nobile},
Theorem 1, Remark 2).


Let us fix some $i \in \{ 1, \dots, m-d \}$. 
If $1 \leq j \leq m$ and $j \notin \mathrm{supp} (n_i)$  then $\ell_j^{(i)}  =0$.
If $1 \leq j \leq m$  and if  $j \in \mathrm{supp} (n_i)$ then 
\[U_j   \cdot \frac{{\partial} h_i}{\partial
U_j} = \left\{ 
\begin{array}{lcl}
\ell_j^{(i)}  U^{n_i^+} & \mbox{ if } & j \in \mathrm{supp} ( n_i^+),
\\
\ell_j^{(i)}  U^{n_i^-} & \mbox{ if } & j \in \mathrm{supp} ( n_i^-).
\end{array}
\right.
\] 
Since $\psi(n_i^+) = \psi (n_i^-)$  we get
$\psi_* ( U^{n_i^+} ) = \psi_*  ( U^{n_i^+} )$.  It follows that the matrix $\tilde{J}$ is the result of multiplying
the $i$th row of the matrix of relations $R$ defined in \eqref{rel2}
by the monomial ${X}^{\psi(n_i^+)}$, 
for $i=1,
\dots p$.

If $1 \leq
j_1 < \cdots < j_d \leq m$, then we get
\begin{equation} \label{clave}
{X}^{\a_{1} + \cdots+  \a_{m}  } \, | {J}_{j_1, \dots, j_d} |=
{X}^{\psi({n_1^+})+ \cdots + \psi ({n_{m-d}^+ } )} {X}^{\a_{j_1}  +
  \cdots+  \a_{j_d}  } \,  |R_{j_1, \dots, j_d}|.
\end{equation}
Formula (\ref{clave}) implies that the Jacobian ideal of
$Z^\Lambda$, which is generated by $ \{ |{J}_{j_1, \dots, j_d} | \}
_{1\leq j_1 < \dots < j_d \leq m}$, and the ideal generated by 
\[ 
\{ {X}^{\a_{j_1}  +
  \cdots+  \a_{j_d}  } \, | R_{j_1, \dots, j_d} | \}_{ 1 \leq j_1 <
\cdots < j_d \leq m }
\]
 are related by invertible ideals,  hence they
have isomorphic blow ups. Finally, by  Lemma \ref{linear}, the determinant $| R_{j_1, \dots,
  j_d}| $ vanishes if and only if $\a_{j_1} \y \cdots  \y \a_{j_d} =0$ , for
$1 \leq j_1 < \cdots < j_d \leq m$.

In the general case, we can assume that 
the  submatrix  $R^{m-d+1, \dots, m}$ is of maximal rank $m-d$.
The Nash modification of $Z^\Lambda$ is 
isomorphic to the blow up of the
ideal $({J}_{j_1, \dots, j_d}^{m-d+1, \dots, m})_{j_1, \dots, j_d}$ (see \cite[Prop. 60]{T3}).
Notice that if we
replace   $R_{j_1, \dots, j_d}$ by
$R_{j_1, \dots, j_d} ^{m-d+1, \dots, m}$, 
and  ${J}_{j_1, \dots, j_d}$ by ${J}_{j_1, \dots, j_d}^{m-d+1, \dots, m}$,  in equation (\ref{clave})  the
resulting formula holds. Then, the proof follows by the previous argument using Lemma \ref{linear}.
\end{proof}

\begin{Defi}
The \textit{logarithmic Jacobian ideal} of $Z^\Lambda$ is the monomial ideal  $\mathcal{J}_\mathrm{log} (Z^\Lambda)$ defined by \eqref{jota-lambda}.
\end{Defi}

\begin{Rem}  
A prime characteristic version of the logarithmic Jacobian ideal of a toric variety was defined in  the paper \cite{MR4792759}. It is also a monomial ideal whose blowing up coincides with the Nash modification. 
\end{Rem}


As a consequence of the proof of Proposition \ref{nash-toric} we get the following description of the generators of 
the Jacobian ideal of the toric variety $Z^\Lambda$.

\begin{Pro}     \label{Jacobian-ideal-toric}
If   $1 \leq j_1 < \cdots  < j_d \leq m $  
 and  $1 \leq i_1 < \cdots  < i_{m-d}  \leq p$
we set
\[
m_{j_1,\dots, j_d}^{i_1,\dots, i_{m-d}} :=
\psi (n_{i_1}^+ ) + \cdots + \psi (n_{i_{m-d}}^+)  - \sum_{j \in \{1, \dots, m\} \setminus \{ j_1, \dots, j_d \}}
 \a_j.
\]
The Jacobian   ideal $\textrm{Jac} (Z^\Lambda)$ of the toric
variety $Z^\Lambda$ is the monomial ideal of $\C [\Lambda]$
generated by 
\[
 ( \{X^{m_{j_1,\dots, j_d}^{i_1,\dots, i_{m-d} }} \mid 
 \a_{j_1} \y \cdots \y \a_{j_d} \ne 0 
 \mbox{ and }   \mathrm{rk} (R^{1, \dots, \hat{i}_1, \dots, \hat{i}_{m-d}, \dots, p} ) = m-d  \} ) 
 \subset  \C [ X^\Lambda].
\]
\end{Pro}
\begin{proof}
Notice that the Jacobian ideal of $Z^\Lambda$ is generated by
the set of minors 
        \[ \{ |{J}_{j_1, \dots, j_d}^{1, \dots, \hat{i}_1,
\dots, \hat{i}_{m-d}, \dots, p} | \} ^{  1 \leq i_1 < \cdots  <  i_{m-d}  \leq p}
_{1\leq j_1 < \dots < j_d \leq m}.\]
 By Formula (\ref{clave}), 
we obtain the equality
\[
 |{J}_{j_1, \dots, j_d}^{1, \dots, \hat{i}_1,
\dots, \hat{i}_{m-d}, \dots, p} | = X^{m_{j_1,\dots, j_d}^{i_1,\dots, i_{m-d} }}
 |R_{j_1,
\dots, j_d} ^{1, \dots, \hat{i}_1, \dots, \hat{i}_{m-d}, \dots, p}
|.
\]
We deduce that $|{J}_{j_1, \dots, j_d}^{1, \dots, \hat{i}_1,
\dots, \hat{i}_{m-d}, \dots, p} | =0$ 
if and only  if $   |R_{j_1,
\dots, j_d} ^{1, \dots, \hat{i}_1, \dots, \hat{i}_{m-d}, \dots, p}
| =0$. Notice that if the rank  
$\mathrm{rk} (R^{1, \dots, \hat{i}_1, \dots, \hat{i}_{m-d}, \dots, p} )$ is less than $m-d$ 
then all the determinants $|R_{j_1,
\dots, j_d} ^{1, \dots, \hat{i}_1, \dots, \hat{i}_{m-d}, \dots, p}|$ vanish.
By Lemma \ref{linear} we have 
$|R_{j_1,
\dots, j_d} ^{1, \dots, \hat{i}_1, \dots, \hat{i}_{m-d}, \dots, p}
| \ne 0$  if and only if 
$\mathrm{rk} (R^{1, \dots, \hat{i}_1, \dots, \hat{i}_{m-d}, \dots, p} ) = m-d$
and $ \a_{j_1} \y \cdots \y \a_{j_d} \ne 0$. 
\end{proof}

\begin{Exam}
Let $\a_1 = (2,0)$, $\a_2 = (0,2)$, $\a_3 = (3,0)$ and $\a_4 = (7,1)$.  One can check that 
$\mathcal{L} = \mathrm{ker}(\psi)$  is a lattice with basis $n_1 = (3,0,-2,0)$ and $n_2 = (7,1, 0,-2)$.
By Proposition \ref{prop: mon eq} below,
$Z^{\Lambda}$ is a complete intersection defined by the equations $h_1 = h_2 = 0$,
where  $h_1 = U_1^3 - U_3^2$ and $h_2 = U_1^7 U_2 - U_4 ^2$. 
We have 
\[
( U_j   \cdot \frac{{\partial} h_i}{\partial U_j})_{i, j} = 
\begin{pmatrix}
3 U_1^3 & 0 & -2 U_3^2 & 0
\\
7 U_1^7 U_2 & U_1^7 U_2 & 0 & -2U_4^2  
\end{pmatrix}
\mbox{ and }
R = 
\begin{pmatrix}
3 & 0 & -2& 0
\\
7  &  1 & 0 & -2 
\end{pmatrix}.
\]
By Proposition \ref{Jacobian-ideal-toric}, the monomial ideal 
of $\C[\Lambda]$ generated by the functions
$X^{m_{1,2}}, X^{m_{1,4}}, X^{m_{2,4}}, X^{m_{3,4}}$ is the Jacobian ideal 
$\textrm{Jac} (Z^\Lambda)$  of $Z^\Lambda$, where 
\[
m_{1,2} = \a_3 + \a_4, 
\quad
m_{1,4} = 7 \a_1 + \a_3, 
\quad 
m_{2,4} = 6 \a_1+ \a_2 + \a_3, 
\quad 
m_{3,4} = 9 \a_1.
\]
\end{Exam}

\section{Quasi-ordinary hypersurface singularities} \label{qo}

A complex analytic germ $(S, 0)$ is {\it quasi-ordinary} if there exists a finite projection
$(S, 0) \rightarrow (\C^d,0)$ which is a local isomorphism outside
a normal crossing divisor.
In the hypersurface case, there is an embedding  $(S, 0) \subset
(\C^{d+1}, 0) $   defined by an equation $f= 0$, where $f
\in \C \{ X \} [X_{d+1}]$ is a {\it quasi-ordinary polynomial}: a
Weierstrass polynomial with discriminant $\D_Y f$ of the form
$\D_Y f = X^\d \epsilon$ for a unit $\epsilon$ in the ring $ \C \{
X \}$ of convergent power series in the variables $X=
(X_1, \dots, X_d)$ and $\d \in \Z^d_{\geq 0}$.


We suppose in this paper that the germ $(S,0)$ is analytically irreducible.
The Jung-Abhyankar
theorem guarantees that the roots of a quasi-ordinary polynomial $f$,
called {\it quasi-ordinary branches},
are fractional power series in the ring $\C \{ X^{1/n} \} $, for
$n=\deg f$ (see \cite{Abhyankar}). If  $\{ \z ^{(l)}
\}_{l =1}^{N} \subset \C \{ X^{1/n} \}$ is the set of roots of $f$,
the discriminant $\D_Y f$ of $f$ with respect to $Y$  is equal to
$
\D_Y f = \prod_{i\ne j} (\z^{(i)} - \z^{(j)}),
$
hence each factor $\z^{(t)} - \z^{(r)}$ is  of the form
$X^{\ll_{t,r}}  \epsilon _{t,r} $ where $\epsilon_{t,r} $ is a unit
in $\C \{ X^{1/n} \} $. The monomials $X^{\ll_{t,r}}$ (resp. the
exponents ${\ll_{t,r}}$) are called {\it characteristic}.
By  Lipman \cite{Lipman2}, the characteristic exponents can be relabelled in the form
\begin{equation} \label{ch-order}
\ll_1 \prec \ll_2 \prec \cdots \prec \ll_g  ,
\end{equation}
where the relation $\prec$ is defined in \eqref{eq-order-def}. The characteristic exponents determine the nested
sequence of lattices
\[
M_0 \subset M_1 \subset \cdots \subset M_g =: M,
\]
where $M_0
:=\Z^{d} $ and  $M_j := M_{j-1} + \Z \ll_j$ for $j=1, \dots, g$.
We set $\ll_0 =0$ and  $\ll_{g +1} = + \infty$.
    The index $n_j$ of the lattice $M_{j-1}$ in the lattice $M_j$ is strictly bigger than one
(see \cite{Lipman2, GP1}).


Lipman noticed that the support of a quasi-ordinary branch has very special properties. 
\begin{Lem} (see \cite[Prop. 1.3]{Gau})  \label{lem:support}
Let $\zeta = \sum c_{\a} X^{\a} \in \C\{ X^{1/n} \}$ be a quasi-ordinary branch with characteristic exponents $\ll_1, \dots, \ll_g$. 
Then, if $c_\a \ne 0$ there exists a unique  $j \in \{ 1, \dots, g+ 1 \}$ such that 
$\ll_{j-1} \preceq \a \npreceq \ll_j$ and then $\a \in M_j$.
\end{Lem}


We say that the quasi-ordinary branch $\zeta$ has {\em well-ordered variables} if the $g$-tuples 
$(\lambda_{1,i}, \dots, \lambda_{g,i})$, which are defined in terms of the coordinates of $\ll_1, \dots, \ll_g$ with respect to the canonical basis, 
satisfy that 
\[ (\lambda_{1,i}, \dots, \lambda_{g,i}) >_{\mathrm{lex}}   (\lambda_{1,j}, \dots, \lambda_{g,j}), \]
for $1 \leq i <  j \leq d$, where $>_{\mathrm{lex}}$ denotes the lexicographic order.
We can relabel the variables $X_1, \dots, X_d$ in order to have this condition.
\begin{Defi} \label{def:normalized}
The quasi-ordinary branch $\zeta$ is {\em normalized} if it is given with well-ordered variables, 
and if $\ll_1$ is of the form 
$\ll_1 = (\ll_{1,1}, 0, \dots, 0)$ then $\ll_{1,1} > 1$. 
\end{Defi}

Lipman proved that the q.o. hypersurface $(S,0)$ can be parametrized by a normalized 
quasi-ordinary branch (see \cite{Gau}). 
The embedded topological type
of $(S, 0) \subset (\C^{d+1}, 0)$ is classified by the characteristic exponents of a 
normalized quasi-ordinary branch $\z$ parametrizing $(S,0)$  (see \cite{Lipman2} and \cite{Gau}).


\subsection{Normalization of a quasi-ordinary singularity}

We study quasi-ordinary hypersurface singularities by using toric geometry
methods.


 \begin{Not} \label{not-can-basis}
 The semigroup $\Z^d_{\geq 0}$ has a minimal
set of generators $e_1, \dots, e_d$, which is a basis $\mathcal{B}$ of the
lattice $M_0$. The dual basis of the dual lattice $N_0$ spans a
regular cone $\s$ in $N_{0,\R}$.    Then,  we have    $\s^\vee
= \R_{\geq 0} e_1 + \dots + \R_{\geq 0} e_d $ and $\Z^d_{\geq 0} =
\s^\vee \cap M_0$. 
If $\a \in \s^\vee$ we write $\a = \sum_{i=1}^d  \a_i e_i$ with respect to the basis $\mathcal{B}$. 
We do this in particular for the characteristic exponents $\ll_j= \sum_{i=1}^d  \ll_{j, i} e_i$.
\end{Not}


The $\C$-algebra of convergent  power series
$\C\{ X_1, \dots, X_d\} $ is isomorphic to $\C \{ \s^\vee \cap M_0 \}$.
This isomorphism identifies the monomial $X_1^{\a_1} \cdots
X_d^{\a_d}$ with the {\em monomial} ${X}^\a  \in  \C [ \s^\vee \cap
M_0]$,     where $\a= \sum_{k=1}^d \a_k e_k$.


We identify the local algebra ${\mathcal O}_S = \C\{X_1, \dots,
    X_d \} [X_{d+1}]/(f)$
of the singularity $(S,0)$ with
the ring  $\C \{X_1, \dots, X_d \} [\z] = \C \{ \s^\vee \cap
    M_0 \} [\z]$.
The normalization of a quasi-ordinary singularity (non necessarily
hypersurface) is analytically isomorphic to a toric simplicial
singularity (see \cite{PP-C}).
In the hypersurface case the normalization is determined by
the characteristic exponents.
\begin{Lem} \label{normal} {\rm (see   \cite[Prop. 14]{GP4})}
The quasi-ordinary branch $\z$ belongs to
$\C \{ \s^\vee \cap M \}$. The homomorphism of $\C$-algebras
${\mathcal O}_S = \C\{ \s^\vee \cap
    M_0\} [\z] \rightarrow \C \{ \s^\vee \cap M \} $
is the inclusion of ${\mathcal O}_S$ in its integral closure ${\mathcal O}_{\overline{S}} = \C \{ \s^\vee \cap M \}$ in
its field of fractions. 
\end{Lem}

\subsection{The semigroup of a quasi-ordinary hypersurface} \label{Sec: semi}

Set
 \begin{equation} \label{rel-semi}
{\overline{\ll}}_{1} =  \ll_1 \quad \mbox{ and } \quad  {\overline{\ll}}_{j+1} = n_j {\overline{\ll}}_{j} +
\ll_{j+1} -  \ll_{j}, \quad \mbox{ for } \quad j= 1, \dots, g-1.
\end{equation}

The following formula can be deduced from \eqref{rel-semi} for $j \in \{2, \dots, g\}$:
\begin{equation} \label{rel-sg-ch}
\overline{\ll}_j = \ll_j + (n_{j-1} -1) \ll_{j-1} + n_{j-1} (n_{j-2} -1) \ll_{j-2}  + \cdots +  n_{j-1} \cdots n_2 (n_{1} -1) \ll_1.
\end{equation}
It implies that 
${\overline{\ll}}_1 \prec \cdots \prec {\overline{\ll}}_g $ by using \eqref{ch-order}.
Formula \eqref{rel-sg-ch} implies also that 
\begin{equation} \label{eq:car-semi}
\ll_{j} \preceq {\overline{\ll}}_{j},\mbox{ for } j = 2, \dots, g.
\end{equation}

\begin{Defi}  (see    \cite{GP1, KM}) \label{barra}
The semigroup
\[ \Gamma := \Z_{\geq 0} e_1 + \cdots + \Z_{\geq 0} e_d + \Z_{\geq 0}  {\overline{\ll}}_1  +  \cdots +  \Z_{\geq 0} {\overline{\ll}}_g 
\subset \s^\vee \cap M \]
is associated to
a sequence of characteristic exponents $\ll_1,\dots, \ll_g$ of a
quasi-ordinary branch $\zeta$. 
By convenience, we denote by     $\g_1,
\dots, \g_{d+g} $ the sequence of generators ${{e}}_1,\dots,
{{e}}_d, {\overline{\ll}}_1,\dots,{\overline{\ll}}_g$ of the semigroup $\Gamma$. We
set also $\g_{g+1} =\infty$.
\end{Defi}

A function $h \in \mathcal{O}_S$ has a  \textit{dominant exponent with respect to 
the quasi-ordinary branch $\zeta$ if 
\[
h =  X^{\g(h)} \cdot u (h), \mbox{ with } \g(h) \in M \mbox{ and } u(h)  \in \C \{ \s^\vee \cap M \}  \mbox{ is a unit.}
\]}
We say that $\g(h)$ is the \textit{dominant exponent} of the function $h$. 
We denote by $\mathcal{C}(\zeta)$ the subset of $\mathcal{O}_S$ consisting of functions 
with a dominant exponent with respect to $\z$.
The set of dominants exponents of functions in $\mathcal{C}(\zeta)$ is equal to the 
semigroup $\Gamma$ (see \cite{PPP04}). The intrinsic nature of the semigroup is more subtle.
\begin{The}  (see \cite{GP1, PPP04,Grenoble,Pseries}).
The semigroup $\Gamma$ is an analytical invariant of the irreducible germ of quasi-ordinary hypersurface $(S,0)$.
\end{The}

\begin{Rem}
If the quasi-ordinary branch $\zeta$ is normalized then 
$e_1, \dots, e_d$, ${\overline{\ll}}_1, \dots, 
{\overline{\ll}}_g$ is a minimal set of generators of $\Gamma$ (see \cite{GP1}). 
Notice that $\Gamma$ is a sub-semigroup of the lattice $M$. In \cite{GP1} we introduced also 
a sub-semigroup ${\overline{\Gamma}}$ of $\Z^{d}$ which is isomorphic to $\Gamma$. 
If $d=1$ the semigroup ${\overline{\Gamma}}$ coincides with the classical semigroup associated to 
the plane branch $(S,0) \subset (\C^2, 0)$. 
Formulas \eqref{rel-semi} and \eqref{rel-sg-ch} are analogous to the classical formulas relating the characteristic of 
a plane branch and the generators of its semigroup (see \cite[Th.3.9]{Zar}).
\end{Rem}

\begin{Rem} \label{rem: app-roots}
One may characterize the generators of $\Gamma$ in terms of 
some of the \textit{approximate roots} of the quasi-ordinary polynomial $f$ (see \cite{GP1}).
Janusz Gwo\'{z}dziewicz and Arkadiusz P\l oski
gave a simplified approach to the Abhyankar-Moh theory of approximate roots
in terms of intersection multiplicities
 (see \cite{AM73, GP95}, see also the survey \cite{PP03}). 
Their presentation was essential in order to extend the results about approximate roots 
to quasi-ordinary hypersurfaces
in the irreducible case in  \cite{GP1}, and more generally in Beata Gryszka's
 paper  \cite{app-Beata}. 
In connection with these results, let us mention 
several irreducibility criterions for quasi-ordinary polynomials
(see \cite{GG21, A12, Villa14, GG15}). 
The study of polar curves was also a relevant topic in
 the work of Arkadiusz P\l oski (see for instance the survey
\cite{GLP}). 
Bunch decompositions of the polars of a quasi-ordinary hypersurface have been given in 
\cite{GBGP05, PP04}. A subtle factorization theorem for higher order polars in the quasi-ordinary hypersurface case was obtained by 
Evelia Garc\'{\i}a~Barroso and Janusz Gwo\'{z}dziewicz in \cite{GBG20}.
\end{Rem}

\begin{Lem} {\rm (see \cite[Lemma 3.3]{GP1})}  \label{poi2}
We have the following properties:
\begin{enumerate} [label=(\alph*)]
\item \label{ref:poi1} $n_{j} \g_{d+j} \prec \g_{d+j+1}$, for $j=2,\dots, g$.

\item \label{poi} We have relations of the form
\begin{equation} \label{can-rel}
 n_j  \g_{d+j}
 =   \ell_{1}^{(j)} \g_1+ \cdots +  \ell_{d}^{(i)} \g_{d}  + \ell_{d+1}^{(j)} \g_{d+1}
 + \cdots + \ell_{d+j-1}^{(j)} \g_{d+j-1}, 
 \end{equation}
where $\ell_{1}^{(j)}, \dots \ell_{d}^{(j)}  \in \Z_{\geq 0}$,
$0 \leq \ell_{d+k}^{(j)} < n_k$ for $k=1,\dots, j-1$,  and 
$j=1,\dots, g$ (see Definition \ref{barra}).
\end{enumerate}
\end{Lem}


As a consequence of  Proposition \ref{prop: mon eq}, the semigroup $\Gamma$ associated with a quasi-ordinary hypersurface is a complete intersection semigroup. This implies that $\Gamma$ has a minimal Frobenius vector $\g_0$, which was studied by Abdallah Assi in \cite{Assi}
(see Section \ref{sec-tor} for the definition of the Frobenius vector).
\begin{Pro} (see \cite{Assi}) \label{def: Frobenius}
The minimal Frobenius vector of the semigroup $\Gamma$ is equal to 
\begin{equation} \label{eq:Frob}
\g_0 :=     \sum_{j=1}^g (n_j -1) \overline{\ll}_{j} - \sum_{i=1}^{d} {e}_j.
\end{equation}
\end{Pro}

\begin{Rem} The \textit{Poincar\'e series}   $P_{\Gamma}(X)  : = \sum_{\g \in \Gamma}  X^{\g} \in \C[[ \Gamma ]] $ of the semigroup $\Gamma$
has a rational form 
\[
P_{\Gamma} (X)  : =
\prod_{j=1}^{g} (1 - X^{n_j \overline{\ll}_{j}}) \prod_{j=1}^{g} (1 -
X^{\overline{\ll}_j})^{-1} \prod_{i=1}^{d} (1 -
X^{{e}_j})^{-1}.
\] 
(see \cite{Pseries}). We observe the following
symmetry property in terms of the  Frobenius vector $\g_0$,
\[
P_{\Gamma}(X)  = (-1)^{d} X^{\g_0}   P_{-\Gamma} (X).
\]
Remark that by Proposition \ref{complete-intersection}, the
variety  $Z^\Gamma$ is a complete intersection, hence $\C[\Gamma]$
is a $M$-graded Cohen Macaulay Gorenstein domain (cf. with
\cite{Stanley} Chapter 1, Theorem 12.7). 
\end{Rem}

\begin{Rem} Poincar\'e series associated with complete intersection semigroups are described similarly in 
\cite[\S 4]{AGO15}, where they are called \textit{Hilbert series}. We follow here the terminology introduced by  
Gusein-Zade, Delgado and Campillo in \cite{GZDC}. See also \cite{CM24} for a survey on some aspects of   Poincar\'e series of semigroups.
\end{Rem}

\subsection{The associated toric variety $Z^\Gamma$}

We consider the equivariant embedding of the toric variety  $
Z^{\Gamma} \subset \C^{d+g}$  given by  $U_i = {X}^{\g_i }$,
for  $i=1,\dots, d+g$, where
$(U_1, \dots, U_{d+g})$ denotes a system of
coordinates on $\C^{d+g}$.
\begin{Pro}   \label{prop: mon eq}
  {\rm  (see \cite[Prop. 38]{GP4})} \label{complete-intersection}
The relations \eqref{can-rel} define the binomials
\begin{equation} \label{canbin}
\begin{array}{cclcl}
h_j & := & U_1^{\ell_{1}^{(j)}}   \cdots
U_{d+j-1} ^{\ell_{d+j-1}^{(j)}} & - &
U_{d+j}^{n_j}, \textrm{ for } j =1, \dots, g,
 \end{array}
\end{equation}
The ideal of the embedding of the toric variety 
$Z^\Gamma \subset \C^{d+g}$ is generated by 
the binomials $h_1, \dots, h_g$. The toric variety  $Z^\Gamma$ is a complete intersection.
\end{Pro}
\begin{Defi}
The matrix of relations associated with (\ref{canbin}) is
\begin{equation} \label{rel}
\left(
\begin{array}{cccccccc}
\ell_{1}^{(1)} &  \cdots & \ell_{d}^{(1)} &  \ell_{d+1}^{(1)}  & 0 & 0&
\cdots & 0
\\
\ell_{1}^{(2)} &  \cdots & \ell_{d}^{(2)} &  \ell_{d+1}^{(2)}  &
\ell_{d+2}^{(2)} & 0& \cdots & 0
\\
\cdots &  \cdots & \cdots &  \cdots  &  \cdots & \cdots & \cdots &
\cdots
\\
\ell_{1}^{(g)} &  \cdots & \ell_{d}^{(g)} &  \ell_{d+1}^{(g)}  &
\ell_{d+2}^{(g)} & \ell_{d+3}^{(g)}  & \cdots &  \ell_{d+g}^{(g)}
\end{array}
\right),
\end{equation}
where we have set $\ell_{d+j}^{(j)} :=  - n_j$ for $j = 1, \dots, g$.
\end{Defi}

We use some elementary properties related with this matrix.

\begin{Defi} \label{erre4}
For $1 \leq i \leq d$ we denote by ${m}_i$ the smallest integer $1 \leq m  \leq g+1$ such that
$\ell_i^{(m)} \ne 0$, where by convention $m_i  = g+1 $ means that  $\ell_i^{(1)} = \dots
=\ell_i^{(g)} =0$.
\end{Defi}

\begin{Lem}  \label{Lem:support2}  Let us take an integer $1 \leq i \leq d$.
\begin{enumerate}[label=(\alph*)]
\item \label{erre3}
If $w \in    \s^\vee \cap  M_{m_i -1}$, 
then $w_i$ is a nonnegative integer.

\item \label{st-b} 
$m_i = \min ( \{ j \in \{ 1, \dots,  g\}  \mid \overline{\ll}_{j,i} \ne 0 \} \cup \{ g +1 \} )$.

\item  \label{proto} If $1 \leq j \leq g$ and 
 $ \ll_{j,i}  \ne 0$,  then we have 
$n_j {\overline{\ll}}_{j,i} \geq 1$. 

\item \label{st-bb}
If $1 \leq j \leq g$ and  $m_i < j$,   then $\overline{\ll}_{j, i} \geq 1$.

\item \label{st-c} If $2 \leq j \leq g$ and if $0 < \overline{\ll}_{j,i}  < 1$,   then $m_i= j$.

\end{enumerate}
\end{Lem}
\begin{proof}
\ref{erre3} The definition of $m_i$ implies that the lattice $M_{m_i-1}$ is the direct sum of 
the sublattice spanned by $e_1, \dots, \widehat{e_i}, \dots, e_d, \ll_1, \dots, \ll_{m_i -1}$
and the sublattice  $\Z e_i$. 
If $w \in M_{m_i-1}$ belongs to the cone $\s^\vee$
spanned by $e_1, \dots, e_d$ then $w_i \in \Z_{\geq 0}$. 

\ref{st-b}
If $1 \leq j < {m}_i$, then the vector $\overline{\ll}_j \in \s^\vee$ belongs to  the lattice 
$ M_j \subset M_{m_i -1} $ by the hypothesis. By  \ref{erre3}, $\overline{\ll}_{j,i}$ is a nonnegative integer. 
By using the direct sum decomposition in \ref{erre3}, if $\overline{\ll}_{j,i} \ne 0$ then
the coefficient $\ell^{(j)}_i$ appearing in \eqref{can-rel} would be nonzero, contradicting the
assumption $1 \leq j < {m}_i$. 
This implies that   $\overline{\ll}_{j,i} = 0$, thus ${\ll}_{j,i}  =0$ by \eqref{eq:car-semi}.
Notice also that 
$n_{m_i} \overline{\ll}_{m_i} \in M_{m_i -1} $ and by definition of $m_i$ 
the coefficient $\ell^{(m_i)}_i$ is nonzero. This implies that  the $i$th coordinate 
$\overline{\ll}_{m_i, i}$ is nonzero.

\ref{proto} If $j=1$  the property holds  since $n_1 {\overline{\ll}}_1 =
n_1 \ll_1 \in M_0 = \Z^d$, thus $ n_1 {\overline{\ll}}_{1, i} >0$ if  $ \ll_{1, i}  \ne 0$. Suppose that the
result is true for $j-1$. 

- If $\ll_{j-1, i} \ne 0$,  by the
induction hypothesis we get $1 \leq n_{j-1}
{\overline{\ll}}_{j-1, i} $.  We get
$ 1 \leq n_{j-1} {\overline{\ll}}_{j-1, i} \leq
{\overline{\ll}}_{j, i}$ by \ref{ref:poi1} in Lemma  \ref{poi2}.

- Otherwise $\ll_{j-1, i} = 0$. By \eqref{rel-semi},
this implies that
 $\ll_{s,i} = 0$, for $s=1, \dots, j-1$.   It follows that $m_i = j$ and $\overline{\ll}_{j,i}= \ll_{j,i} > 0$ by \eqref{rel-semi}.
Since $n_j \overline{\ll}_j \in M_{j-1}$  by (\ref{can-rel}), we  deduce from 
statement \ref{erre3} that
$n_j \overline{\ll}_{j, i} \geq 1$.

\ref{st-bb}  
By \ref{st-b} we have $\overline{\ll}_{m_i, i} \ne 0$. By \ref{proto} we get 
$n_{m_i}  \overline{\ll}_{m_i, i}  \geq 1$. Since $m_i < j$,  it follows from \ref{st-b} that 
$\overline{\ll}_{j, i} \geq  n_{m_i}  \overline{\ll}_{m_i, i}  \geq 1$ by \eqref{rel-semi}.


\ref{st-c} 
Notice that $\overline{\ll}_j = n_{j-1} \overline{\ll}_{j-1} + \ll_j - \ll_{j-1}$ by 
\eqref{rel-semi}, 
where $ 0 \prec \ll_j - \ll_{j-1}$ by \eqref{ch-order}. If $\overline{\ll}_{j-1, i} \ne 0$ then we get 
$\overline{\ll}_{j,i} \geq n_{j-1} \overline{\ll}_{j-1, i} \geq 1$ by \ref{proto}, thus we must have 
$\overline{\ll}_{j-1, i} = 0$. 
It follows that 
$\overline{\ll}_{1, i} = \cdots = \overline{\ll}_{j-1, i} = 0$. Then, the assertion follows by \ref{st-b}.
\end{proof}

\section{The logarithmic Jacobian ideal of a quasi-ordinary hypersurface}\label{rel-Jacobian}

We introduce the logarithmic Jacobian ideal of a singularity with toric normalization
following  \cite{LR} in the toric case.  We review the normal toric case
following \cite[Chapter 3]{Oda} and \cite[Appendix]{LR}.
Then, we apply this to the case of an irreducible germ of quasi-ordinary hypersurface and 
also to its associated monomial variety.

\subsection{Logarithmic jacobian ideals}
Let $Y$ be the algebroid germ defined by an ideal $I$ of $\C [[X_1, \dots, X_n ]]$ with analytic algebra
$A = \C [[X_1, \dots, X_n ]] / I$.
We denote by  $\Omega_Y^1$ the $A$-module  of {\em K\"ahler differential forms}
and by $d : A \rightarrow  \Omega_Y ^1$ its canonical derivation.
   As usual we denote by $ \Omega^k_Y$ the $A$-module
$ \Omega^k_Y := \bigwedge_k \Omega^1_Y$.      See \cite{GLS} Chapter I, \S 1.10.


We consider the toric singularity $Z=Z^{\s^\vee \cap M}$
 with formal local algebra $\mathcal{O}_{Z}= \C [[\s^\vee \cap M]]$ at the origin.
We denote by $D$ the equivariant Weil divisor defined by the sum of orbit closures of codimension one
in the toric variety $Z$.
The $\mathcal{O}_{Z}$-module   $\Omega_{Z}^1 (\log D) $
 of $1$-forms of $Z$ with logarithmic poles along $D$ is
identified with $ \mathcal{O}_Z \otimes_\Z M $. We have a map of
$\mathcal{O}_Z$-modules 
\[
  \eta: \Omega_{Z}^1 \rightarrow
\mathcal{O}_Z \otimes_\Z M , \quad \mbox{ such that }  d X^{\g} \mapsto
X^\g \otimes \g, \mbox{ for } \g \in \s^\vee \cap M.
\] 
If $u_1, \dots, u_d $ is a basis of the lattice $M$ we write the expansion of $\g \in  \s^\vee \cap M$ 
as $\g = \sum_{i=1}^d \g_i u_i$. 
If $h = \sum_{\g \in \s^\vee \cap M} c_\g X^\g \in  \mathcal{O}_{Z}$ then we have $
\eta (d h ) =  \sum_{\g \in \s^\vee \cap M} c_\g X^\g \otimes  \g
= \sum_{i=1}^d \left( \sum_{\g \in   \s^\vee \cap M}     c_\g \, \g_i X^\g \right)  \otimes u_i$.

The map $\y^k \eta$  is a homomorphism of
$\mathcal{O}_{Z}$-modules such that
\[
\begin{array}{lcl}
\y^k \eta : \Omega_{Z}^k &  \longrightarrow  &    \Omega_{Z}^k (\log D) = \mathcal{O}_Z \otimes_\Z \bigwedge_k  M,
\\
 d X^{\g_1} \y \cdots \y    d X^{\g_k} &  \mapsto &   X^{\g_1 + \cdots + \g_k }
\otimes  \g_1 \y \cdots \y \g_k,
\end{array}
\]
for $k = 1, \dots, d$.
Let us fix an orientation of $M$ and assume that the basis  $u_1,  \dots,  u_d$ is positively oriented. We have an isomorphism $ \bigwedge_d M
\rightarrow \Z$, which sends $u_1 \y \cdots  \y u_d $ to $1$.
This
induces an identification  $\Omega_{Z}^d (\log D) \equiv \mathcal{O}_Z$. It follows that $\eta^d (\Omega_Z^d)$ corresponds
by this identification to an ideal of  $\mathcal{O}_Z$, which is independent
of the basis of $M$ chosen. This ideal is called the {\em
logarithmic Jacobian ideal} of $Z$ in \cite{LR}. It is a monomial
ideal generated by $( X^{\b_{j_1} + \cdots +  \b_{j_d}} )_{
\b_{j_1} \y \cdots \y \b_{j_d} \ne 0 } ^{1\leq j_1 < \cdots < j_d
\leq r } $, where $\{ \b_1, \dots, \b_r \}$ is a generating system of the semigroup $\s^\vee \cap M$.

\begin{Defi} \label{def:lji}
Let  $W$ be an analytically irreducible algebroid germ of singularity of dimension $d$ such
that its normalization
$Z$ is a toric singularity.
The {\em logarithmic Jacobian ideal $\mathcal{J}_{\mathrm{log}} (W) \mathcal{O}_{Z} $}  is
   the  ideal of $\mathcal{O}_{Z}$
      generated by the image of the module of differentials
$\Omega_{W}^d$ by the composite $\t$ of the maps
\begin{equation}              \label{composite}
\Omega_{W}^d \longrightarrow  \Omega_{Z}^d \stackrel{\wedge_d \eta}{\longrightarrow}
\Omega_{Z}^d (\log D) \equiv \mathcal{O}_{Z},
\end{equation}
where   $\Omega_{W}^d \rightarrow  \Omega_{Z}^d$ is the map induced by the normalization.
\end{Defi}

         \subsection{The logarithmic Jacobian ideal of the quasi-ordinary hypersurface $(S,0)$.}

We consider an analytically irreducible quasi-ordinary hypersurface $(S,0)$ germ as in Section \ref{qo} and 
we describe its logarithmic Jacobian ideal using that its normalization is a toric singularitiy.

\begin{Defi}   \label{def:eq-xi}
We set  $\mathcal{J}_g $ the  monomial ideal of 
$\mathcal{O}_{Z}$ generated by monomials $X^{\a}$ for $\a$ in the set 
\begin{equation} \label{eq-xi-j}
\Xi_g := 
  \{ \sum_{j=1}^d e_j   \} \cup \{  e_1 + \cdots + \widehat{e_i} + \cdots +  e_d + \ll_{m_i} 
\mid i \in \{1, \dots, d \},  m_i \in \{ 1, \dots, g \} 
\}. 
\end{equation}
We denote by $\mathcal{I}_g$ the  monomial ideal of 
$\mathcal{O}_Z$ generated by $X^{\a}$ for $\a$ in the set 
\begin{equation} \label{eq-xi-j-0}
  \{ \g_{i_1}  + \cdots + \g_{i_d}  \mid 1 \leq i_1 < \cdots < i_d \leq d+g, \quad \g_{i_1}  \y \cdots \y \g_{i_d} \ne 0 \}, 
\end{equation}
(see Definition \ref{barra}).
\end{Defi}

\begin{The} \label{jotak2} (see \cite[Prop. 9.2]{pseries-qo}) 
The ideal $\mathcal{J}_{\mathrm{log}} (S) \mathcal{O}_Z$
is  equal to the  monomial ideal $\mathcal{J}_g$.
\end{The}

\begin{proof}
We have the inclusion   ${\iota} : \mathcal{O}_{{S}} \rightarrow
\mathcal{O}_{\overline{S}}$ corresponding to the normalization map,
where  $\mathcal{O}_{\overline{S}}  =  \C[[\s^\vee \cap M]]$ by Lemma
\ref{normal}. We consider a fixed quasi-ordinary branch $\z = \sum_\a c_\a X^\a$  parametrizing $(S,0)$.
We denote by $x_i$ the image of the coordinate $X_i$
in $\mathcal{O}_{S}$,  for $i=1, \dots d+1$. By Lemma
\ref{normal}, the function ${\iota}(x_{d+1}) = \z $ belongs to $\C\{ \s^\vee \cap
M \} $ and ${\iota} (x_i) = X^{{e}_i}$ for $i=1, \dots, d$ (see Notation \ref{not-can-basis}).


We analyze the
image of $ d x_{j_1} \y \cdots \y d x_{j_d}$ by the composite $\t$
of the maps  (\ref{composite}) for $W = S$. 
 Since $ d x_{j_1} \y
\cdots \y d x_{j_d}$ generate the   $ \mathcal{O}_{S}$-module
$\Omega_{S}^d$   it follows that $\t    (d x_{j_1} \y \cdots \y
d x_{j_d})$ generate the ideal $\t
(\Omega_{S}^d)$ of $ \C \{ \s^\vee \cap M \} $,  for $1 \leq j_1 < \dots < j_d \leq d+1$.

If $(j_1, \dots, j_d) = (1, \dots, d)$,  then  $\t    (d x_{1} \y \cdots \y d x_{d}) $ is equal to
$X^{e_1 + \dots + e_d}$ times a nonzero constant, hence $X^{e_1 +
\dots + e_d} $ belongs to $\mathcal{J}_{\mathrm{log}} (S)\mathcal{O}_Z$.  

Otherwise, we have $j_d = d+1$ and $ 1 \leq j_1 < \dots <  j_{d-1} \leq d$. 
We denote by $i$ the integer such that $\{ j_1, \dots, j_{d-1}, i \} = \{ 1, \dots, d\}$. 
Let us consider the index $m_i$  introduced in  Definition \ref{erre4}.
We decompose the quasi-ordinary branch $\z = \sum_{\alpha} c_\a X^\a$ as  $\z =  \z_{i}^-
   +  \z_{i}^+$,  where   $ \z_{i}^+ =\sum_{\ll_{m_i} \preceq \a} c_\a  X^\a$  and 
   $\z_{i}^- = \z - \z_{i}^+ $ (see Lemma \ref{lem:support}).

By definition we have $\z_{i}^- =  \sum_{\ll_{m_i} \npreceq \a} c_\a  X^\a$. 
Let $\a \in M$ be an element appearing in the support of the series $\z_{i}^-$.
By Lemma \ref{lem:support},  $\a$ belongs to the lattice $M_{m_i -1}$.
In addition, if we expand $\a= \sum_{s=1}^d \a_s e_s$ in terms of the canonical basis $\mathcal{B}$,  then 
the coefficient $\a_i$ is a nonnegative integer by Lemma \ref{Lem:support2} \ref{erre3}.
We write
\[
\begin{array}{lcl}
\t    (d x_{j_1} \y \cdots \y d x_{j_d}) & = &   \wedge_d \eta  (
d X^{{e}_{j_1}} \y \cdots 
  \y  d X^{{e}_{j_{d-1}}} \y d \z_{i}^{-})
\\
&  &  +   \wedge_d \eta  ( 
d X^{{e}_{j_1}} \y \cdots \y  d X^{{e}_{j_{d-1}}} \y d \z_{i}^{+})
\end{array}
 \]  
Then, the term
\[ \wedge_d \eta  ( d X^{{e}_{j_1}} \y \cdots 
  \y  d X^{{e}_{j_{d-1}}} \y d \z_{i}^{-})  = 
\sum_{\ll_{m_i} \npreceq \a} c_\a  \cdot e_{j_1}\y \cdots 
\y e_{j_{d-1}} \y  \a  \cdot  X^{e_{j_1} + \cdots + e_{j_{d-1}} +  \a }
\]
is divisible by $X^{e_1 +\cdots + e_d}$. 
We deduce that 
\[ \wedge_d \eta  ( d X^{{e}_{j_1}} \y \cdots \y  d X^{{e}_{j_{d-1}}} \y d \z_{i}^{+})  = 
\sum_{\ll_{m_i} \preceq \a} c_\ll  \cdot e_{j_1}\y \cdots 
\y e_{j_{d-1}} \y  \a  \cdot  X^{e_{j_1} + \cdots + e_{j_{d-1}} +  \a } 
\]
is the product of the monomial $X^{e_{j_1} + \cdots + e_{j_{d-1}} + \ll_{m_i}}$ times a unit, since 
\[ c_{\ll_{m_i}}  \cdot e_{j_1}\y \cdots 
\y e_{j_{d-1}} \y  \ll_{m_i} \ne 0.\]
This implies that 
$\mathcal{J}_{\mathrm{log}} (S) \mathcal{O}_Z$
is equal to  the monomial ideal $\mathcal{J}_g $
(see Definition \ref{def:eq-xi}).
\end{proof}

Let $\Lambda \subset M$ a subsemigroup generated by $\a_1, \dots, \a_m \in M$ such that 
$\Z \Lambda = M$ and $\R_{\geq 0} \Lambda = \s^\vee$ as in Section \ref{sec-tor}. 
Then, the canonical map $Z := Z^{\s^\vee \cap M} \rightarrow Z^\Lambda$ is the normalization of 
$Z^\Lambda$.  Notice that the logarithmic Jacobian ideal $\mathcal{J}_{\mathrm{log}} (Z^\Lambda) \mathcal{O}_Z$ is the extension of the ideal of $\mathcal{O}_{Z^\Lambda}$ defined in  Proposition
\ref{nash-toric}. The following result is a
consequence of the definitions and the previous discussion in the toric case. 

\begin{Pro} (see also \cite{LR})    \label{log-tor}
The  ideal $\mathcal{J}_{\mathrm{log}} (Z^\Lambda) \mathcal{O}_Z$ 
is the monomial ideal of $Z^{\s^\vee \cap M}$ defined by the monomials
\[
\{ X^{\a_{j_1}  +
 \cdots+  \a_{j_d}}   \, \mid \, 1 \leq j_1 < \cdots < j_d
 \leq m \mbox{ and }  \a_{j_1}  \y
 \cdots \y \a_{j_d} \ne 0 \}.
\]
\end{Pro}


We apply Proposition \ref{log-tor} to 
the particular case of $Z^\Gamma$ taking  into account Definitions  \ref{barra} and  \ref{def:eq-xi}.
\begin{Cor}
The logarithmic Jacobian  ideal $\mathcal{J}_{\mathrm{log}} (Z^\Gamma) \mathcal{O}_Z$ is equal to the monomial ideal $\mathcal{I}_g$. 
\end{Cor}

\begin{Rem}
Lejeune and Reguera have
shown that the logarithmic Jacobian ideal of a normal affine toric
surface plays a significant role in the rational expression of the
associated {\em geometric motivic Poincar\'e series} (see  \cite{LR}). In a joint work with Helena
Cobo, we have shown similar results in the cases
of affine toric varieties \cite{pseries-toric} and of
irreducible quasi-ordinary hypersurface singularities \cite{pseries-qo}. 
\end{Rem}

\section{Description of the Jacobian ideal} \label{Sec:JI}

In this section we give some relations between the Jacobian ideal of the quasi-ordinary hypersurface $S$, and the logarithmic Jacobian ideals of $S$ and of its associated monomial variety $Z^\Gamma$. We need to introduce first a deformation of $S$ which has generic fiber isomorphic to $S$ and specializes to $Z^\Gamma$.

\subsection{Overweight  deformations}   \label{eq-def}

In \cite{Tesis} we construct an embedding \[ (S, 0) \hookrightarrow
(\C^{d+g},0)\] together with a deformation
$\mathcal{S}$ of $S$ over $Z^\Gamma$, with generic fiber $S_p \subset \C^{d+g}$ isomorphic to $S$ and
special fiber $S_0 = Z^\Gamma \subset \C^{d+g}$. This deformation is equisingular in the sense that there exists
a toric modification $\p: W \rightarrow \C^{d+g}$,
characterized by the properties of the semigroup $\Gamma$, such that
it provides a simultaneous normalization of  $(S_p, 0)$, and  the composite of $\p$ with any toric resolution of the singularities of $W$
is an embedded resolution of $(S_p, 0)$, for $p=0$ or $p$ in the torus of $Z^\Gamma$ (see \cite{GP4}).
This result generalizes one of
Teissier and Goldin for plane branches (see \cite{Goldin}).
See also
Teissier's development of this approach in terms of valuations
\cite{VDT, MR3329046}.


Let us describe this deformation. 
First, we have the exact sequence \eqref{eq:sequence} associated to the map
\[
\psi:\Z^{n+g} \to M, \quad \psi(u_i) =  \g_i, \mbox{ for } i = 1, \dots, d+g,
\]
which sends the elements of the canonical basis $u_1, \dots, u_{d+g}$ of 
$\Z^{n+g}$ to the generators of the semigroup $\Gamma$ (see \eqref{eq:presentation} and Definition \ref{barra}).

We consider the map $\psi^*: \C[U_1, \dots, U_{d+g}] \to \C[\Gamma]$  as in \eqref{iota}.
We define the {\em weight} of a monomial $U^\alpha \in  \C[U_1, \dots, U_{d+g}] $ by $\psi ({\alpha}) \in \Gamma$. This defines a $\Gamma$-grading of $\C[U_1, \dots, U_{d+g}]$. 
Notice that the binomials $h_j$ defined in \eqref{canbin} 
are weighted homogeneous of degree $n_j \overline{\ll}_j$, for $j=1, \dots, g$.


\begin{Defi} An overweight deformation of the binomials $h_j$  in \eqref{canbin}, for $j = 1, \dots, g$, is of the form 
\begin{equation} \label{cantri}
H_{j} = h_j    + c_j \cdot t^{\overline{\ll}_{j+1} - n_j \overline{\ll}_j}  \cdot
U_{d+j+1} + \sum_{n_j \overline{\ll}_j \prec \psi (\alpha^{(j)})}    c_{\alpha^{(j)}} \cdot t^{\psi(\a^{(j)}) - n_j\overline{\ll}_j} \cdot
U_1^{ \a^{(j)}_1} \cdots  U_{d+j}^{\a^{(j)}_{d+j}},  
\end{equation}
where ${\alpha^{(j)}} = \sum_{i=1}^{d+i} {\alpha^{(j)}_i} u_i$ and the coefficients verify 
$c_{\alpha^{(j)}} \in \C$, $c_j \in \C^*$ for $j=1, \dots,g-1$, and $c_g=0$.  
Notice that  $H_{j}$ is an element of the ring $\C[t^\Gamma][[U_1, \dots, U_{d+g}]]$ for $j =1, \dots, g$.
\end{Defi}
The terms added to $h_j$ in \eqref{cantri} have weight higher than $n_j \overline{\ll}_j$. That is why we say that  $H_j$ is an overweight deformation of $h_j$. We use Teissier's terminology in \cite{MR3329046}.


If $p$ is a closed point in $Z^\Gamma$ we can evaluate the monomials of the form $t^\g$, for 
$\g \in \Gamma$, appearing in the expansion \eqref{cantri}. We obtain the series
\begin{equation} \label{cantri-p}
H_{j,p} = h_j    + c_j \cdot t^{\overline{\ll}_{j+1} - n_j \overline{\ll}_j} (p)   \cdot
U_{d+j+1} + \sum_{\a^{(j)}}    c_{\alpha^{(j)}} \cdot t^{\psi(\a^{(j)}) - n_j \overline{\ll}_j}  (p) \cdot
U_1^{ \a^{(j)}_1} \cdots  U_{d+i}^{\a^{(j)}_{d+i}}. 
\end{equation}
Note that  $H_{j,p}  \in \C[[U_1, \dots, U_{d+g}]]$.
The ideal $(H_{1,p}, \dots, H_{g,p})$ defines a germ $S_p$ at the origin of $\C^{d+g}$, which may be seen as the fiber of a family $\mathcal{S}$ defined by 
 $(H_{1}, \dots, H_{g})$ over the point 
 $p \in Z^\Gamma
$.
Next, we describe some of fibers $S_p$ 
of such a family, when $p$ is the origin ${0}$
of the toric variety $Z^{\Gamma}$ or when $p$ belongs to the torus of $Z^\Gamma$.

\begin{Pro} \label{emb} (see \cite[Prop. 3.15, 3.17]{Tesis} and \cite[Prop. 39]{GP4}).
Let $(S,0)$ be a germ of quasi-ordinary singularity with associated semigroup $\Gamma$. 
There exists an overweight deformation  $H_1, \dots, H_g$ of  $h_1, \dots, h_g$ such that:
\begin{enumerate}[label=(\alph*)]
\item If $p$ belongs to the torus of $Z^\Gamma$ 
then the  germ $(S_p,0)$ is analytically isomorphic to $(S,0)$. 
\item If $p =0$ then $(S_0, 0) = (Z^\Gamma, 0) \subset (\C^{d+g},0)$.
\end{enumerate}
\end{Pro}

If $p$ belongs to the torus of $Z^\Gamma$, or if $p = 0 \in Z^\Gamma$, we have
a map
\begin{equation}  \label{bar-iota}
{\iota}_p : \C \{ U_1,\dots, U_{d+g} \}  \rightarrow \C \{ \s^\vee \cap M\},
\end{equation}  
defined as the composition of the
canonical maps   $\C\{ U_1,\dots, U_{d+g} \} \rightarrow \mathcal{O}_{S_p} $ with
the inclusion of  $\mathcal{O}_{S_p}$ in its integral closure $ \mathcal{O}_{\overline{S}_p}$, and 
then
\begin{equation} \label{eq-int-clos}
 \mathcal{O}_{\overline{S}_p} =  \C\{ \s^\vee \cap M\}.
\end{equation}
 For $p$ in the torus the equality \eqref{eq-int-clos} is consequence of Lemma \ref{normal} and Proposition \ref{emb}.
 See \cite{GP4} for the case $p=0$ of the equality \eqref{eq-int-clos}.

\begin{Rem} \label{rem:iota}
Take an overweight deformation  $H_1, \dots, H_g$ of  $h_1, \dots, h_g$ and a point $p $ in the torus of $Z^\Gamma$. Then
 \[ t^{\overline{\ll}_{j+1} - n_j \overline{\ll}_j} (p) \ne 0, \mbox{ for } j =1, \dots g-1.\]
 We can use this fact to eliminate recursively the variables 
 $U_{d+g}, \dots, U_{d+2}$ in the equations
(\ref{cantri}). We set $X_i := U_i$, for $i =1, \dots, d+1$ and we obtain the equation of the embedding $(S_p, 0) \subset (\C^{d+1}, 0)$, with coordinates $(X_1, \dots, X_{d+1})$. By \cite{GP4} the restriction to $S_p$ of the projection 
\[ (X_1, \dots X_{d+1}) \rightarrow (X_1, \dots, X_d) \] is a quasi-ordinary projection, 
and 
$S_p$ is parametrized by a quasi-ordinary branch $X_{d+1} = \z_p$, which has characteristic monomials 
$\ll_1, \dots, \ll_g$.
 It follows that  $\mathcal{J}_{\mathrm{log}} (S_p) = \mathcal{J}_g $, for all $p $ in the torus of $Z^\Gamma$ (see Definition \ref{def:eq-xi}).
In particular, if $\underline{1}$ is the unit point of the torus $Z^{M}$, we can set
 $(S,0) := (S_{\underline{1}}, 0)$, and then  the statement of Proposition \ref{emb} holds for this deformation. 
 \end{Rem}

\subsection{The extension of the matrix of relations}

We consider an overweight deformation $H_1, \dots, H_g$ of the binomials $h_1, \dots h_g$ 
defining the toric variety $Z^\Gamma$. It defines  a deformation  $\mathcal{S}$,
with fibers $S_p$, for $p \in Z^\Gamma$.
\begin{Not} \label{not: jmatrix}
We denote by $J(p)$ (resp. $\tilde{J}(p)$) the matrix with coefficients
in the ring $\C\{\s^\vee \cap M \}$, which is defined from the  Jacobian matrix of $(H_{1,p}, \dots,
H_{g, p})$ by the formulas
\[
J(p) := \left( {\iota} (\frac{{\partial} H_{j,p}}{\partial U_i})
\right)_{i=1, \dots, d+g}^{j=1, \dots, g}, \quad \left( \textrm{resp. }
\tilde{J}(p) := \left( {\iota} (U_i \frac{{\partial}
H_{j,p}}{\partial U_i}) \right)_{i=1, \dots, d+g}^{j=1, \dots, g} \right),
\]
where ${\iota}$ is the map defined by (\ref{bar-iota}).
We denote by $R(p)$ the matrix obtained from $\tilde{J}(p)$
by factoring out from the entries of the  $j$th row,  the term
       $X^{   n_j \overline{\ll}_j} $, for $j=1, \dots, g$.
We denote the entry $L^{(j)}_i  (p) $ of the matrix $R(p)$ in the row $j$ and in the column $i$, 
simply by $L^{(j)}_i $
to avoid cumbersome notation. 
\end{Not}

\begin{Lem} \label{lem:mext}
The matrix $R (0)$ has integer coefficients and coincides with the matrix of relations 
\eqref{rel} associated with the presentation of the semigroup $\Gamma$. 
\end{Lem}
\begin{proof}
By Proposition \ref{emb} the terms appearing in $H_j$ have
weight greater than or equal to $   
n_j \overline{\ll}_j$, with equality precisely for those terms appearing in $h_j$.
This implies that 
 $X^{   n_j \overline{\ll}_j} $ divides the entries of the $j$th row of
$\tilde{J}(p)$, for $j=1,   \dots, g$.
It follows that the matrix 
$R(p) $ is of the form
 \begin{equation} \label{eq:Rp}
\left(
\begin{array}{cccccccc}
L_{1}^{(1)} &  \cdots 
&  L_{d+1}^{(1)}  &
\epsilon_{1} {X}^{{\overline{\ll}}_2 -
  n_{1} {\overline{\ll}}_{1}}   & 0 & 0 &  0
\\
L_{1}^{(2)} &  \cdots 
&  L_{d+1}^{(2)}  &
L_{d+2}^{(2)} & \epsilon_{2}  {X}^{{\overline{\ll}}_3 -
  n_{2} {\overline{\ll}}_{2}}  &  0 & 0
\\
\cdots &  \cdots 
 &  \cdots  &  \cdots
& \cdots & \cdots &  \cdots
\\
L_{1}^{(g-1)} &  \cdots 
&  L_{d+1}^{(g-1)}  &
L_{d+2}^{(g-1)} & L_{d+3}^{(g-1)}  & \cdots 
& \epsilon_{g-1} {X}^{{\overline{\ll}}_g -
  n_{g-1} {\overline{\ll}}_{g-1}}
\\
L_{1}^{(g)} &  \cdots 
&  L_{d+1}^{(g)}  &  L_{d+2}^{(g)}
& L_{d+3}^{(g)}  & \cdots  
&  L_{d+g}^{(g)}
\end{array}
\right),
\end{equation}
where $L_i^{(j)} \in \C \{ \s^\vee \cap M \}$ and  the coefficient 
\begin{equation} \label{eq:epsilon}
\epsilon_{j}  : = c_j \cdot t^{\overline{\ll}_{j+1} - n_j \overline{\ll}_j} (p) 
\end{equation}
vanish when $p= 0$, while it is nonzero if $p$ belongs to the torus of $Z^\Gamma$.

The
 coefficient $L_{i}^{(j)}$  has an expansion of 
 the form
\begin{equation} \label{elle}
L_{i}^{(j)} = \ell_{i} ^{(j)} + \sum_{n_j \overline{\ll}_j \prec \psi (\alpha^{(j)}), \, \alpha^{(j)}_i \geq 1} 
c_{\alpha^{(j)}} \cdot t^{\psi(\a^{(j)}) - n_j \overline{\ll}_j}  (p) \cdot
X^{\psi (\a^{(j)}) -  n_j {\overline{\ll}}_j},
\end{equation}
for $j = 1, \dots,g$ and $i =1,\dots, d+j$. 
Notice that the condition $\alpha^{(j)}_i \geq 1$ characterizes when the partial derivative
$\frac{\partial}{\partial U_i}  (U_1^{ \a^{(j)}_1} \cdots  U_{d+i}^{\a^{(j)}_{d+i}})$  is nonzero.
The weight condition  $n_j \overline{\ll}_j \prec \psi (\alpha^{(j)})$ implies 
that $t^{\psi(\a^{(j)}) - n_j \overline{\ll}_j}  (0) = 0$. 
\end{proof}

\begin{Defi}
We say that the matrix $R (p)$ is the {\em extension of the matrix of relations} $R(0)$.
\end{Defi}

The following lemma is a consequence of the proof of Lemma \ref{lem:mext}.
\begin{Lem}\label{elle2}
If  $1 \leq j  < {m}_i$ the monomial
${X}^{e_i}$ divides $ L_{i}^{(j)}$,  for $i=1, \dots, d$ (see Notation
\ref{erre4}).
\end{Lem}
\begin{proof}
Consider the expansion 
(\ref{elle})  of $L_i^{(j)}$. 
Take $\a^{(j)} \in \Z^{d+g}_{\geq 0}$ such that 
$n_j \overline{\ll}_j \prec \psi (\alpha^{(j)}) $  and  $\alpha^{(j)}_i \geq 1$.
By Lemma \ref{Lem:support2} \ref{st-b}, we have  $\overline{\ll}_{j,i} =0$.
Then, the condition $\alpha^{(j)}_i \geq 1$
implies that $e_i \prec \psi (\alpha^{(j)})  - n_j \overline{\ll}_j $, that is, 
$X^{e_i}$ divides $X^{\psi (\alpha^{(j)})  - n_j \overline{\ll}_j }$ in the ring 
$\C\{ \s^\vee \cap M \}$.
This implies the assertion.   
\end{proof}

\subsection{Description of the image of the Jacobian ideal of $S_p$}

Next, we provide a description of the Jacobian ideal of $S_p$ in terms of the matrix $R(p)$.
The approach  is inspired by Proposition \ref{Jacobian-ideal-toric} in the toric case.

\begin{Pro} \label{prop:jac}
The extension  $\mathrm{Jac}(S_p) \mathcal{O}_Z$ of the Jacobian ideal of $S_p$ in the local algebra of the normalization $Z$ of $S_p$ is  generated by 
\[
 \{ 
 X^{\g_0 + \g_{i_1} + \cdots + \g_{i_d}} \cdot
 |R(p)_{i_1, \dots, i_d} |
 \mid 
 1 \leq i_1 < \cdots < i_d \leq d+g \}. \]
 \end{Pro}
\begin{proof}
The Jacobian ideal of $S_p$ can be expressed in terms of the generators 
$H_{1, p}, \dots, H_{g, p}$ of the defining ideal of $S_p$. 
It is the ideal of the local algebra 
$\mathcal{O}_{S_p}$ defined by the 
 vanishing of the minors of order $g$ of the Jacobian matrix $(\frac{\partial H_{j,p}}{\partial U_i})$.
It follows that the ideal $\mathrm{Jac}(S_p) \mathcal{O}_Z$  is defined by the minors 
$|J(p)_{i_1, \dots, i_d}|$ 
of order $g$ of the matrix $J(p)$, for $1 \leq i_1 < \cdots < i_d \leq d+g$ (see Notation \ref{not: jmatrix}). 
Then, we argue as in the proof of \eqref{clave} to obtain the formula
\begin{equation} \label{clave-tilde}
{X}^{\g_1  + \cdots+  \g_{d+g}  } \,  | J_{i_1, \dots, i_d} (p) |=
{X}^{n_1 {\overline{\ll}}_1 + \cdots + n_g {\overline{\ll}}_g } {X}^{\g_{i_1}  +
  \cdots+  \g_{i_d}  } \,  |{R(p)}_{i_1, \dots, i_d} |.
\end{equation}
The assertion follows by 
taking into account the description of the Frobenius vector $\g_0$ (see Proposition \ref{def: Frobenius}). 
\end{proof}

\begin{Rem} \label{rem:pzero}
By Lemma \ref{linear}, if $p = 0 $ is the origin of $ Z^\Gamma$ we have
\begin{equation} \label{eq: nonzero}
|{R(0)}_{i_1, \dots, i_d} | = 0 \mbox{ if and only if }  \g_{i_1}  \y \cdots 
\y  \g_{i_d}  \ne 0.
\end{equation} 
 We obtain the following formula, which connects 
the Jacobian ideal and the logarithmic Jacobian ideal
\[
\mathrm{Jac}(S_0) \mathcal{O}_Z = X^{\g_0} \cdot \mathcal{J}_{\mathrm{log}} (S_0)  \mathcal{O}_Z =  X^{\g_0} \cdot  \mathcal{I}_g.
\]
It is a particular case of Proposition  \ref{Jacobian-ideal-toric} (see Remark \ref{log-tor} and Definition \ref{def:eq-xi}). 
\end{Rem}

\begin{Rem} \label{rem: multiplicity}
Let us apply the previous discussions to the case of plane curve singularities, which was one of the favourite
research topics of Arkadiusz P\l oski.
We consider the case when $S$ is a plane branch, that is, when $d=1$.
Assume that the line $X_1 = 0$ is not tangent to $S$. 
 The algebra of the normalization of $S$ is of the form $\mathcal{O}_Z = \C \{ \s^\vee \cap M  \}$,  where
$ \s^\vee \cap M = \frac{1}{n} \Z_{\geq 0}$ and $n$ is the multiplicity of $S$. 
This implies that
$\g_1, \g_2, \dots, \g_{g+1} \subset \frac{1}{n} {\Z_{\geq 0}}  $ 
is the minimal system of generators of the semigroup $\Gamma$ (see Definition \ref{barra}).

Since $d =1$, $\mathcal{O}_Z$ is a discrete valuation ring, and if 
 $t := X^{\frac{1}{n} }$ then  $\mathcal{O}_Z = \C \{ t \}$.
We denote by $\nu_t$ its $t$-adic valuation.
Let us denote by  $\overline{\b}_0, \dots, \overline{\b}_g \in \N$ 
the minimal system of generators 
of the \textit{classical semigroup} of the branch $S$, defined by 
taking intersection multiplicities of $S$ with other curves not containing
the branch $S$ as a component (see \cite{Zar} or \cite{Pf-Pl}). 
We have $\nu_t ( X^{\g_i}) = n \g_i = \overline{\b}_{i-1}$, for $i =1, \dots, g+1$.
The {\em multiplicity  of the Jacobian ideal} of the branch $S$ is 
defined by 
\[
e( \mathrm{Jac} (S) ) = \nu_t ( \mathrm{Jac} (S) \mathcal{O}_Z).
\]
Let us explain how our method provides a formula for $e( \mathrm{Jac} (S_p) )$ for $p \in Z^\Gamma$.


In this case, the matrix $R(p)$ is of order $g \times (g+1)$. By Proposition \ref{prop:jac}, the ideal
$\mathrm{Jac} (S_p) \mathcal{O}_Z$ is generated by 
$
\{ 
 X^{\g_0 + \g_{i}}  \cdot
 |R(p)_{i} | \mid i = 1, \dots g+1 \}$.
By definition we can expand the minor $|R(p)_{i} |$ as a series in  $  \C\{ t \}$ with constant term $|R(0)_{i} |$.
By  \eqref{eq: nonzero} the minor $|R(0)_{i} |$ does not vanish since $\g_i \ne 0$, for $i =1, \dots, g+1$.  
Notice also that 
$ \min_{i=1, \dots, g+1} \{ \g_0 + \g_i \} = \g_0 + \g_1$.
It follows that $\mathrm{Jac} (S_p) \mathcal{O}_Z = (  X^{\g_0 + \g_{1}} ) \mathcal{O}_Z$. 
We deduce the following equality
\begin{equation} \label{eq:jac-mult}
e( \mathrm{Jac} (S_p) ) = \nu_t (X^{\g_0 + \g_{1}}) =  n (\g_0 + \g_1) \stackrel{\eqref{eq:Frob}}{=} \sum_{j=1}^g (n_j -1) \overline{\b}_{j}.
\end{equation}
This formula holds in particular when $p = \underline{1}$ is the unit point of the torus of $Z^\Gamma$, and when $p = 0$. These points define the fibers $S_{\underline 1} = S$ and $S_0 = Z^\Gamma$.
A particular case of a result of Teissier in \cite[Prop. II.1.2]{T73} in the case of plane branches provides the relation 
\begin{equation} \label{eq:milnor}
\mu(S) = e( \mathrm{Jac} (S) ) - n +1, 
\end{equation}
where $\mu(S)$ is the \textit{Milnor number} and $n$ is the multiplicity of the branch $S$. See  \cite{LeDT}
 and \cite{Greuel}  for the case of complete intersection curves. Formula \eqref{eq:jac-mult} can be deduced from 
equation \eqref{eq:milnor} by using the expression of the Milnor number of a branch in terms of 
the minimal system of generators of its semigroup \cite[Chapter II.3]{Zar}, see also \cite[Section 15]{GBPPGP25}. Let us mention the contributions
of Arkadiusz  P\l oski and Evelia Garc\'\i a Barroso on the Milnor number of curves in arbitrary characteristic
(see for instance \cite{MR3839793}).
\end{Rem}

We come back to the general case of a quasi-ordinary hypersurface $S$ of dimension $d$.

\begin{Pro} \label{th:inclusion1}
For any point $p$ in the torus of $Z^\Gamma$ we have
\[
X^{\g_0} (\mathcal{I}_{g} + \mathcal{J}_g ) \subset \mathrm{Jac}(S_p) \mathcal{O}_Z. 
\]
\end{Pro}
\begin{proof}
It is enough to prove the inclusions $X^{\g_0} \cdot \mathcal{I}_{g} \subset \mathrm{Jac}(S_p) \mathcal{O}_Z$ and 
$X^{\g_0} \cdot \mathcal{J}_{g} \mathcal{O}_Z  \subset \mathrm{Jac}(S_p) \mathcal{O}_Z$.


Let us prove the inclusion $X^{\g_0} \cdot \mathcal{I}_{g} \subset \mathrm{Jac}(S_p) \mathcal{O}_Z$ first.
Take $1 \leq i_1 < \cdots < i_d \leq d+g$  such that $\g_{i_1}  \y \cdots 
\y  \g_{i_d}  \ne 0$. We obtain a generator
$X^{\g_{i_1} + \cdots \g_{i_d}}$
 of the ideal $\mathcal{I}_g$ (see Definition \ref{def:eq-xi}).  Then, we have
 \[
 X^{\g_0 + \g_{i_1} + \cdots + \g_{i_d}} \cdot
 |R(p)_{i_1, \dots, i_d} |  \in 
 \mathrm{Jac}(S_p) \mathcal{O}_Z
 \] by Proposition \ref{prop:jac}.
Notice that 
$|R(p)_{i_1, \dots, i_d} |$ is a unit of $\mathcal{O}_Z = \C[[\s^\vee \cap M]]$, since 
it can be expanded as a power series with constant term 
equal to $|R(0)_{i_1, \dots, i_d} | $, which is nonzero by the hypothesis and \eqref{eq: nonzero}.
Thus, the monomial $X^{\g_0 + \g_{i_1} + \cdots + \g_{i_d}} $ belongs to $\mathrm{Jac}(S_p) \mathcal{O}_Z$. 
 

We prove the inclusion $X^{\g_0} \cdot \mathcal{J}_g \subset \mathrm{Jac}(S_p) \mathcal{O}_Z$. 
Notice that by the previous case $X^{\g_0 + \g_1 + \cdots \g_d}$ belongs to  $\mathrm{Jac}(S_p) \mathcal{O}_Z$. 
By Definition \ref{def:eq-xi}, it is enough to prove that 
for any  vector $\a$ of the form 
\[
\a =e_1 + \cdots + \widehat{e_i} + \cdots +  e_d + \ll_{m_i} 
\mbox{ , where } i \in \{ 1, \dots, d \} \mbox{ and } m_i \in \{1, \dots, g\},
\]
 the monomial 
$X^{\g_0 + \a} $ belongs to  $\mathrm{Jac}(S_p) \mathcal{O}_Z$.


If $m_i = 1$ then we have 
 $e_1 \y  \cdots \y \widehat{e_i} \y \cdots \y e_d \y \ll_{1} \ne 0$. Since $\ll_1 = \overline{\ll}_1$ it follows that 
  $X^{e_1 + \cdots + \widehat{e_i} + \cdots +  e_d + \ll_{1}}  $ belongs to 
  $\mathcal{I}_g$
and the assertion follows by the previous case (see Definition  \ref{barra}).


We assume that  $1 < m_i \leq g$ and we set 
\begin{equation} \label{eq:delta}
\d:= \g_0 + \g_{1} + \cdots + \widehat{\g_{i}} +   \cdots + \g_{d} + \g_{d+1}.
\end{equation}
We use that 
\begin{equation} \label{eq:term-jac}
 X^{\d}
 |R(p)_{1, \dots, \hat{i}, \dots, d+1} | \in \mathrm{Jac}(S_p) \mathcal{O}_Z,
\end{equation}
by Proposition \ref{prop:jac}.
We expand the minor $ |{R(p)}_{1, \dots, \hat{i}, \dots, d+1 } |$ by its first column and we get 
\[
|{R(p)}_{1, \dots, \hat{i}, \dots, d+1 } | = \sum_{r=1}^g (-1)^{r+1} L_i^{(r)} |{R(p)}_{1, \dots, d+1 }^r |.
\]

If $1 \leq r < m_i$,  then $X^{e_i}$ divides $L_i^{(r)}$  by Lemma \ref{elle2}. It follows that the term
$X^{\d}   L_i^{(r)} |{R(p)}_{1, \dots, d+1 }^r | $ is divisible by the monomial $X^{\g_0 + \g_1 + \dots + \g_d} \in \mathrm{Jac}(S_p) \mathcal{O}_Z$, thus 
\begin{equation}  \label{eq:term-jac2}
\sum_{i=1}^{m_i-1} X^{\d}   L_i^{(r)} |{R(p)}_{1, \dots, d+1 }^r |  \in \mathrm{Jac}(S_p) \mathcal{O}_Z.
\end{equation}


If $m_i \leq  r \leq g $,  then we use the structure of the matrix $R(p)$ (see \eqref{eq:Rp}).
The matrix ${R(p)}_{1, \dots, d+1 }^r $ is block lower triangular. It has  
 two blocks in the diagonal
\[
A_r (p)  := R(p)^{r, \dots, g}_{1, \dots, d+1, d+ r+1, \dots, d+g}, \quad \mbox{ and} \quad
B_r (p)  =  R(p)^{1, \dots, r}_{1,\dots,d+r}, 
\]
where 
\[
A_r (p)  = \begin{pmatrix}
\epsilon_{1} X^{\overline{\ll}_{2} - n_1 \overline{\ll}_{j_1}} & 0 & \dots & 0 
\\
\dots & \dots & \dots    & \dots
\\
* & * & \dots &  \epsilon_{r-1} X^{\overline{\ll}_{r} - n_{r-1} \overline{\ll}_{r-1}}
\end{pmatrix}
,
\]
and 
\begin{equation} \label{eq-br0}
B_r (0)  =  \begin{pmatrix}
n_{r+1} & 0 & \dots & 0 
\\
\dots & \dots & \dots    & \dots
\\
* & * & \dots & n_g 
\end{pmatrix}
.
\end{equation}

Notice that the coefficients \eqref{eq:epsilon} are nonzero if $p$ belongs to the torus of $Z^\Gamma$. Then, the equality 
\[
 \overline{\ll}_2 - n_1 \overline{\ll}_1 + \cdots +  \overline{\ll}_{r} - n_{r-1} \overline{\ll}_{r-1}  \stackrel{\eqref{rel-semi}}{=}  
{\ll_2} - {\ll_1} + \cdots +  {\ll_{r}} - {\ll_{r-1}} = \ll_r - \ll_1, 
\]
implies that the determinant $|A_r (p)|$ is of the form $X^{\ll_r - \ll_1}$ times a unit.

The determinant $|B_r (p)|$ is a power series in the ring $\C[[\s^\vee \cap M]]$ which has nonzero constant term $|B_r (0)| = n_{r+1} \cdots n_g$ by \eqref{eq-br0}.
By using \eqref{ch-order}, there exists a unit $\epsilon \in \C[[\s^\vee \cap M]]$ such that 
\[ 
\sum_{r=m_i}^g (-1)^{r+1} L_i^{(r)} |{R(p)}_{1, \dots, d+1 }^r |
= X^{\ll_{m_i} - \ll_1} \cdot \epsilon.
\]
By \eqref{eq:term-jac} and \eqref{eq:term-jac2} we deduce that 
$
 X^{\d + \ll_{m_i} - \ll_1} \cdot \epsilon \in \mathrm{Jac}(S_p) \mathcal{O}_Z
 $.  
 This implies the assertion,  since by \eqref{eq:delta} and Definition \ref{barra} we have
 \[ 
 \begin{array}{lcl}
 \d + \ll_{m_i} - \ll_1 & = &  \g_0 + \g_{1} + \cdots + \widehat{\g_{i}} +   \cdots + \g_{d} + \g_{d+1} + \ll_{m_i} - \ll_1 
 \\ 
 & =  & 
 \g_0 + e_{1} + \cdots + \widehat{e_{i}} +   \cdots + e_{d} +  \ll_{m_i}.
 \end{array}
 \] 
 \end{proof}

Next, we prove the following inclusion of the monomial ideals in Definition  \ref{def:eq-xi}.
\begin{Pro}  \label{prop:mon-inc}
We have $\mathcal{I}_g \subset \mathcal{J}_g$ 
\end{Pro}

\begin{proof}
Let us take $1 \leq j_1 < \cdots < j_k \leq d$ and $1 \leq i_{k+1} < \cdots < i_{d} \leq d$ such that 
\begin{equation} \label{eq: wedge}
\overline{\ll}_{j_1} \y \cdots \y \overline{\ll}_{j_k} \y e_{i_{k+1}} \y \cdots \y e_{i_{d}} \ne 0. 
\end{equation}
Set $\a = \sum_{\ell=1}^k \overline{\ll}_{j_\ell} + \sum_{\ell=k+1}^d e_{i_{\ell}}$.
It is enough to prove that for any $\a$ of this form there exists an element $\b \in \Xi_g$  (see  \eqref{eq-xi-j}) such that $\b \preceq \a$.  For instance, if $k=0$ then $\a = \sum_{\ell =1}^d e_\ell \in \Xi_g$ and 
we take $\b = \a$.
Up to relabelling,  we may assume that $i_{\ell} = \ell$, for $\ell \in \{ k+1, \dots, d \}$ and then
\[ 
 \a = \sum_{\ell=1}^k \overline{\ll}_{j_\ell} + \sum_{\ell=k+1}^d e_{{\ell}}, \] 
(we do not assume here that the quasi-ordinary branch $\zeta$ has well-ordered variables).
If $k=1$ then 
$\a = \overline{\ll}_{j_1} + e_2 + \cdots + e_d $, and $\ll_{j_1, 1} \ne 0$ by \eqref{eq: wedge}. 
By Lemma \ref{Lem:support2} \ref{st-b}, we have $m_1 \leq j_1$. Then, we can take
$\b := \overline{\ll}_{m_1} + e_2 + \cdots + e_d \in \Xi_g$, 
 and  $\b \preceq \a$.


Next, assume  that $k > 1$.  Condition \eqref{eq: wedge} means that the  matrix \[
C := (\overline{\ll}_{j_\ell, i})_{i=1, \dots, k}^{\ell= 1, \dots, k}. 
\] 
has rank $k$.


\textbf{Claim}.   If $i \in  \{ 1, \dots, k \}$ then 
\begin{equation} \label{eq: ine-mi}
m_i \leq j_k,
\end{equation}
and equality in \eqref{eq: ine-mi} happens for at most one index $i_0 \in  \{ 1, \dots, k \}$.


Notice  the inequality $m_i >  j_k$ implies that 
the $i$th column of the matrix $C$ is zero, and this contradicts the condition 
$|C| \ne 0$   (see 
Lemma \ref{Lem:support2} \ref{st-b}). 
If $m_{i_0} = j_k$, the $i$th column vector of the matrix $C$ is of the form 
$(0, \dots, 0 , \overline{\ll}_{j_k, i_0})$ with $\overline{\ll}_{j_k, i_0} \ne 0$. 
Notice that we cannot have two columns of this form 
since $|C| \ne 0$. This ends the proof of the claim.


$\bullet$ If the inequality \eqref{eq: ine-mi} is strict for all  $i \in  \{ 1, \dots, k \}$, then
by Lemma \ref{Lem:support2} \ref{st-bb}, the inequality 
\begin{equation} \label{eq: key-1}
\overline{\ll}_{j_k, i} \geq 1, 
\end{equation}
holds for $i = 1, \dots, k$. 
Then, $\b = e_1 + \cdots + e_d $ belongs to $\Xi_g $ and  $\b \preceq \a$. This ends the proof in this case. 

$\bullet$ If the equality in \eqref{eq: ine-mi} holds for one index $i_0$, we suppose up to relabelling that $i_0 =k$. By Lemma \ref{Lem:support2} \ref{st-bb}, the inequality \eqref{eq: key-1} holds for $i = 1, \dots, k-1$. 


We distinguish two subcases.


- If $\overline{\ll}_{j_k, k} \geq 1$, then we can take $\b = e_1 + \cdots + e_d \in \Xi_g$, 
and $\b \preceq \a$ as in the previous case.


- Otherwise, we have  $0 < \overline{\ll}_{j_k, k}  < 1$. Then, we can take 
$\b = \ll_{j_k} + e_1 + \cdots + e_{k-1} + e_{k+1} + \cdots + e_d \in \Xi_g$,  since 
$m_k = j_k$.
We 
claim that  
$\b \preceq \a$. 
It is equivalent to prove that 
\begin{equation} \label{eq: keysim}
\ll_{j_k} + \sum_{\ell= 1}^{k-1} e_\ell \preceq \sum_{\ell=1}^k \overline{\ll}_{j_\ell}. 
\end{equation}
Notice that the $k$th column vector  of the matrix $C$ is 
$(0, \dots, 0 , \overline{\ll}_{j_k, k})$.  This implies that the first minor of order $k-1$ of $C$ is nonzero.
As in the previous case, we use this fact to prove that 
if $i \in \{ 1, \dots, k-1 \}$ then 
\begin{equation} \label{eq: ine-mi-2}
m_i \leq j_{k-1},
\end{equation}
with equality for at most one index $i_1$. 
We distinguish then two subcases:


- If the inequality in \eqref{eq: ine-mi-2} is strict  for all $i\in \{ 1, \dots, k-1\}$ 
then 
$\overline{\ll}_{j_{k-1}, i} \geq 1$  by Lemma \ref{Lem:support2} \ref{st-bb}.
This implies that 
$ e_1 + \cdots + e_{k-1} \preceq \overline{\ll}_{j_{k-1}}$ and 
\eqref{eq: keysim}
 holds since $\ll_{j_k} \preceq \overline{\ll}_{j_k}$. 


-  Otherwise, we have an equality in  \eqref{eq: ine-mi-2}  for 
one index $i_1$.  We may suppose that $i_1 = k-1$. 
If $\overline{\ll}_{j_{k-1}, k-1} \geq 1$ then \eqref{eq: keysim} holds 
since \eqref{eq: ine-mi-2} implies that 
 $\overline{\ll}_{j_{k-1}, i} \geq 1$
for $i \in \{1, \dots, k-2 \}$ by Lemma \ref{Lem:support2} \ref{st-bb}.


It remains to proof the case $0 < \overline{\ll}_{j_{k-1}, k-1} <1$. 
We set $j = j_{k}$ and $p = j_{k-1}$ for simplicity.  We compute using  \eqref{rel-sg-ch}
\[ 
\begin{array}{lcl}
\overline{\ll}_j  - \ll_j + \overline{\ll}_p & = & \displaystyle{\sum_{\ell =1}^{p -1}}  (n_{\ell} -1) n_{\ell +1} \cdots 
 n_{p-1} \cdot ( 1+ n_{p} \cdots n_{j-1}) \ll_\ell
 \\
 & + & \left( 1 + (n_{p} -1) n_{p +1} \cdots n_{j-1} \right) \ll_{p} 
\\
 & + &   \displaystyle{\sum_{\ell = p +1}^{j -1}} (n_\ell -1)  n_{\ell +1} \cdots n_{j-1} \ll_\ell. 
\end{array}
\]
Notice that the coefficient 
$a_p :=  1 + (n_{p} -1) n_{p +1} \cdots n_{j-1}$ of $\ll_p$ in this formula verifies that
$
a_ p  > n_p 
$, 
since $a_p  - n_p = (n_{p +1} \cdots n_{j-1}-1) (n_p -1)  \geq 1$. 

Take $i \in \{ 1, \dots, k-1 \}$. 
Then we have 
\begin{equation} \label{eq:llji}
\overline{\ll}_{j,i}  - \ll_{j,i} + \overline{\ll}_{p,i} \geq a_p \,  {\ll}_{p,i} > n_p \, {\ll}_{p,i} = 
n_{j_{k-1}} \ll_{j_{k-1}, i},
\end{equation}
where the last equality is just obtained by rewriting. 


If $i \in \{ 1, \dots, k-2 \}$,  then we have a strict inequality in 
\eqref{eq: ine-mi-2}, and by Lemma \ref{Lem:support2} \ref{st-bb} we get  $\overline{\ll}_{j_{k-1}, i} \geq 1$. 
In addition, by
 Lemma \ref{Lem:support2} \ref{proto}, the condition $0 < \overline{\ll}_{j_{k-1}, k-1} <1$ implies $\overline{\ll}_{j_{k-1}, k-1}  = {\ll}_{j_{k-1}, k-1} $ and 
$
n_{j_{k-1}} \ll_{j_{k-1}, k-1} \geq 1$.
By combining these inequalities with \eqref{eq:llji} we get 
\begin{equation} \label{eq: keysim2}
\overline{\ll}_{j_k,i}  - \ll_{j_k,i} + \overline{\ll}_{j_{k-1} ,i} \geq 1, 	 \mbox{ for } i \in \{1, \dots, k-1 \}.
\end{equation}
This means that 
\[ \ll_{j_k}  + e_1 + \cdots + e_{k-1}  \preceq \overline{\ll}_{j_k}  + \overline{\ll}_{j_{k-1}} \preceq 
\overline{\ll}_{j_k}  + \overline{\ll}_{j_{k-1}} + \dots + \overline{\ll}_{j_{1}} .\] 
Therefore,  \eqref{eq: keysim} holds and this ends the proof.
\end{proof}

 As a corollary of Theorem \ref{jotak2},  Propositions \ref {th:inclusion1} and  \ref{prop:mon-inc} we deduce the following result.
\begin{The} \label{th:inclusion}
For any point $p$ in the torus of $Z^\Gamma$ we have 
 \[
X^{\g_0} \, 
\mathcal{J}_{\mathrm{log}} (S_p)  \subset \mathrm{Jac}(S_p) \mathcal{O}_Z. 
\]
\end{The}

 Next, we study the case of quasi-ordinary surfaces.
 
\begin{Pro} \label{prop:jac2}
Assume that $d =2$ and  that $S$ is parametrized by a normalized quasi-ordinary branch $\zeta$ (see Definition \ref{def:normalized}).
Take a point $p$ in the torus of $Z^\Gamma$. 
Let $I$ be the ideal  of $\mathcal{O}_Z$ 
generated by 
\[
 \{ 
 X^{\g_{i_1} + \g_{i_2}} \cdot
 |R(p)_{i_1, i_2} |
 \mid 
1 \leq i_1 < \cdots < i_2 \leq 2+g \} .
 \]
Then, we have the inclusion
$I \subset  \mathcal{J}_g$.
 \end{Pro}
\begin{proof} 
Since $\zeta$ is normalized, we have
\begin{equation} \label{eq:normalized}
\ll_{1,1} \ne 0 \mbox{ and if } \ll_{1,2} = 0 \mbox{ then } \ll_{1,1} > 1.
\end{equation} 
If  $1 \leq i_1 < i_2 \leq 2+ g$ and if  $\g_{i_1} \y \g_{i_2} \ne 0$ then we have 
$X^{\g_{i_1} + \g_{i_2}} \cdot
 |R(p)_{i_1, i_2} | \in \mathcal{I}_g$. 
 We use here that $\mathcal{I}_g \subset \mathcal{J}_g $ by Proposition \ref{prop:mon-inc} 
 

 By the hypothesis \eqref{eq:normalized} and  \eqref{ch-order}, 
if $\g_{i_1} = \overline{\ll}_{j_1}$ and 
$\g_{i_2} = e_2$ then we get $\overline{\ll}_{j_1} \y e_2 \ne 0$. 
Thus, if $\g_{i_1} \y \g_{i_2} =  0$ then there are two cases:

\begin{enumerate}[label=(\alph*)]
\item \label{case 1}
Case  $i_1 = 1$ and $i_2 = 2 + j_1$ with $1 \leq j_1 \leq g$ and  $e_1 \y \overline{\ll}_{j_1} = 0$.
\item \label{case 2}
Case $i_1 = 2 + j_1$, $i_2 = 2 + j_2$ with $1 \leq j_1 < j_2 \leq g$ and  $\overline{\ll}_{j_1} \y
 \overline{\ll}_{j_2} = 0$.
\end{enumerate}


\noindent
\textbf{Case \ref{case 1}}. The case $j_1 = 1$ was considered in the proof of Proposition \ref{th:inclusion1}.
It is enough to consider the case $j_1 > 1 $.  Recall the definition of
 the index $m_2$ introduced in Definition \ref{erre4}. 
By the hypothesis $e_1 \y \overline{\ll}_1 = 0$, we get $e_1 \y {\ll}_1 = 0$, thus
$j_1 < m_2$. We expand the minor 
$ |R(p)_{1, 2+j_1} |$ by its first column and we get 
\[
|R(p)_{1, 2 + j_1} | = \sum_{j=1}^{g} L_2^{(j)} |R(p)_{1,2, 2 + j_1}^{j}|.
\]

We distinguish two subcases: 

- If $1 \leq j < m_2$ we have  $X^{e_2}$ divides the coefficient $L_2^{(j)}$
 (see Lemma \ref{elle2}). It follows that 
$X^{e_1 + e_2}$ divides the term $X^{\g_{1} + \g_{2+j_1}} L_2^{(j)} |R(p)_{1,2, 2 + j_1}^{j}|$, hence 
this term belongs to $\mathcal{I}_g$. This ends the proof of this case if $m_2 = g+1$.


\noindent
- If $ m_2 \leq j \leq g$, the matrix $|R(p)_{1, 2 , 2+j_1}^{j}|$ is block lower triangular. It has three blocks in the diagonal
$A_j$, $B_j$ and $C_j$ where
\[
A_j = R(p)_{1, 2, 2 +j_1, \dots, 2+g}^{{j_1}, \dots, g}, 
\quad 
B_j = R(p)_{1,2, \dots, 2 + j_1, 2+j+1, \dots, 2+g }^{1, \dots, j_1 -1, j, \dots,g},
\quad 
C_j = R(p)_{1, \dots,2+j}^{1, \dots, j}.
\]
In particular, the matrix $B_j$ is diagonal of the form 
\begin{equation} \label{eq:Bj}
B_j = 
\begin{pmatrix}
\epsilon_{j_1} X^{\overline{\ll}_{j_1 +1} - n_{j_1} \overline{\ll}_{j_1}} & 0 & \dots & 0 
\\
\dots & \dots & \dots    & \dots
\\
* & * & \dots &  \epsilon_{j-1} X^{\overline{\ll}_{j} - n_{j-1} \overline{\ll}_{j -1}}
\end{pmatrix}
.
\end{equation}
Notice that $\sum_{s=j_1}^{j-1} \overline{\ll}_{s +1} - n_{s} \overline{\ll}_s =  \sum_{s=j_1}^{j-1}
 {\ll}_{s +1} -  {\ll}_s = \ll_{j} - \ll_{j_1}$.
We deduce that $X^{\ll_j -\ll_{j_1}}$ divides the minor $|R(p)_{1, 2 , 2+j_1}^{j}|$. 
It follows that $X^{e_1 + \overline{\ll}_{j_1}} |R(p)_{1, 2 , 2+j_1}^{j}| $ is divisible by $X^{e_1 + \ll_j}$, since $\ll_{j_1} \preceq \overline{\ll}_{j_1}$. Since $m_2 \leq j $ we have  $e_1 \y \ll_j \ne 0$, hence the monomial $X^{e_1 + \ll_j}$ belongs to $\mathcal{J}_g$. This ends the proof in this case.
 
 
\noindent
 \textbf{Case \ref{case 2}}  
 By hypothesis $\ll_{j_1, 1} \ne 0$, and since $j_1 < j_2$,  
we get  $\overline{\ll}_{j_2,1} \geq 1$.

Assume first that  $\ll_{j_1, 2} \ne 0$. By Lemma \ref{Lem:support2} \ref{proto} we get $n_{j_1} 
\overline{\ll}_{j_1, 2} \geq 1$. Since $j_2 > j_1$, we obtain that 
$\overline{\ll}_{j_2,2} \geq \overline{\ll}_{j_1, 2} \geq 1$ by definition \eqref{rel-semi} and relation \eqref{ch-order}.
We have shown that 
$e_1 + e_2 \prec \overline{\ll}_{j_2}$,  therefore
 $X^{\overline{\ll}_{j_1} + \overline{\ll}_{j_2}} |R(p)_{2+j_1, 2 + j_2} |$ belongs to $\mathcal{I}_g$, since it is divisible 
by $X^{e_1+e_2} \in \mathcal{I}_g$.


 Otherwise we must have $\ll_{j_1, 2} = 0$ and then  $\ll_{j_2, 2} = 0$, since $\overline{\ll}_{j_1} \y \overline{\ll}_{j_2} = 0$, and by definition of $m_2$ we have also $j_1 < j_2 < m_2$. 
 We expand the minor $|R(p)_{2+j_1, 2 + j_2} |$ by its second column and we get 
 \[
|R(p)_{2+ j_1, 2 + j_2} | = \sum_{j=1}^{g} L_2^{(j)} |R(p)_{2,2+ j_1, 2 + j_2}^{j}|.
\]
We distinguish two subcases:

\noindent
-  If $1 \leq j < m_2$ we have  $\ll_{j,2} = 0$ by Lemma \ref{Lem:support2} \ref{st-b}.  Then, 
the monomial $X^{e_2}$ divides the coefficient $L_2^{(j)}$ (see Lemma \ref{elle2}).
Since $\ll_{j_1,1} \ne 0$, we get $\overline{\ll}_{j_2,1} \geq 1$ thus $e_1 \prec \overline{\ll}_{j_2}$.
The monomial $X^{e_1+ e_2} \in \mathcal{I}_g$ divides term $X^{\overline{\ll}_{j_1} + \overline{\ll}_{j_2}} L_2^{(j)} |R(p)_{2,2+ j_1, 2 + j_2}^{j}|$, hence  this term belongs to $\mathcal{I}_g$. This ends the proof when 
$m_2 = g+1$. 

\noindent 
- If $m_2 \leq j \leq  g $, 
the matrix $R(p)_{2,2+ j_1, 2 + j_2}^{j}$ is block lower triangular. It has 
 three blocks in the diagonal $A_j$, $B_j$ and $C_j$ where
\[
A_j = R(p)_{2, 2 +j_1, 2+j_2, \dots, 2+g}^{{j_2}, \dots, g}, 
\quad 
B_j = R(p)_{1, \dots, 2 + j_2, 2+j+1, \dots, 2+g }^{1, \dots, j_2 -1, j, \dots,g},
\quad 
C_j = R(p)_{1, \dots,2+j}^{1, \dots, j}.
\]
In particular, $B_j$ is diagonal of the form  \eqref{eq:Bj} replacing $j_1$ by $j_2$.
Therefore,  the minor $|B_j|$ is equal to the monomial $ X^{\ll_{j} - \ll_{j_{2}}}$ times a unit. As in the previous case, the term $X^{\overline{\ll}_{j_1} + \overline{\ll}_{j_2}} |R(p)_{2,2+ j_1, 2 + j_2}^{j}|$ is divisible by $X^{e_1 + \ll_j}$, since $\ll_{j_2} \preceq \overline{\ll}_{j_2}$. Since $m_2 \leq j $ we have $e_1 \y \ll_j \ne 0$, hence the monomial $X^{e_1 + \ll_j}$ belongs to $\mathcal{J}_g$. This ends the proof in this case.
\end{proof}
 

The following Corollary is an immediate consequence of Propositions \ref{prop:jac}, \ref{prop:jac2},  and Theorems \ref{jotak2},   \ref{th:inclusion}.
 \begin{Cor} \label{cor:inclusion1}
 Assume that the dimension of the quasi-ordinary hypersurface $S$ is equal to two.
 The inclusion 
\[ \mathrm{Jac}(S_p) \mathcal{O}_Z = X^{\g_0} \mathcal{J}_{\mathrm{log}} (S_p),\]
holds for any point $p$ in the torus of $Z^\Gamma$.
\end{Cor}

\begin{Rem}
It is an open question to determine if  the inclusion in 
Theorem 
\ref{th:inclusion}
is an equality when the  dimension of $S$ is greater than $2$. 
\end{Rem}

\begin{Cor}    \label{nor-nash}
Let $S$ be an analytically irreducible quasi-ordinary hypersurface singularity of dimension two. 
The composite
of the normalization map of $S$ with the normalized blow up of the monomial ideal 
 $\mathcal{J}_{\mathrm{log}} (S)$
 is the normalized Nash modification of
$S$.
 \end{Cor}
 \begin{proof}
 Since $S $ is a complete intersection, the Nash
modification of $S$ is isomorphic to the blow up of the Jacobian
ideal in $\mathcal {O}_{S}$ (see \cite{Nobile}). Then the
normalized Nash modification of $S$ is equal to the composition
of the normalization map with the normalized blow up of
${\textrm{Jac}(S)} \mathcal{O}_{Z} $
(see \cite{LT}, Propositions 3.2 and 3.3). Then, the statement
follows from Corollary \ref{cor:inclusion1},
since the ideals
 ${\textrm{Jac}(S)} \mathcal{O}_{Z} $  and $
\mathcal{J}_{\mathrm{log}} (S)  \mathcal{O}_{Z} $ are related by invertible monomial ideals, hence 
they have isomorphic blow ups.
 \end{proof}

\begin{Exam} \label{ej:dos}
Let $(S, 0)$ be a quasi-ordinary hypersurface parametrized by a quasi-ordinary branch  $\z$ with characteristic exponents
\[
\ll_1 = (\frac{3}{2}, 0), \quad \ll_2 = (\frac{7}{4}, 0), \quad \ll_3 = ( 2, \frac{1}{2}).
\]
The characteristic integers are $n_1 = n_2 = n_3 = 2$. 
The generators of the semigroup $\Gamma$ are
\[
e_1 = (1, 0) \quad, e_2 = (0,1), \quad \overline{\ll}_1 =  (\frac{3}{2}, 0), \quad  \overline{\ll}_2 = (\frac{13}{4}, 0),  \quad  \overline{\ll}_3= (\frac{27}{4}, \frac{1}{2}).
\]
Notice that the group $M$ generated by $\Gamma$ is equal to the lattice with basis $v_1 := \frac{1}{4} e_1$, $v_2 :=\frac{1}{2} e_2$. 
The saturation of $\Gamma$ in $M$ is the semigroup $ \s^\vee \cap M = \Z_{\geq 0} v_1 + \Z_{\geq 0} v_2$. 
This implies that the normalization $Z = Z^{\s^\vee \cap M}$  of $S$ is smooth, and $\mathcal{O}_Z $ is isomorphic to 
the ring of convergent power series $\C\{ z_1, z_2 \}$ 
where $z_1 = X^{v_1}$ and $z_2 =X^{v_2}$.
For simplicity, we write the coordinates of the generators of the semigroup $\Gamma$ (resp. of the characteristic exponents) with respect to the basis $v_1, v_2$.
We obtain
\[
\begin{array}{c}
e_1 = (4, 0) \quad, e_2 = (0,2), \quad \overline{\ll}_1 =  (6, 0), \quad  \overline{\ll}_2 = ({13}, 0),  \quad  \overline{\ll}_3= ({27},{1}) 
\\
\mbox{ (resp. }  {\ll}_1 =  (6, 0), \quad \ll_2 = (7,0), \quad \ll_3 = (8,1)).
\end{array}
\]
The minimal Frobenius vector with respect of this basis is $\g_0 = (42, -1)$.
The following relations hold
\[
2 \overline{\ll}_1 = 3 e_1, \quad 2 \overline{\ll}_2 = 5 e_1 + \overline{\ll}_1, \quad 2 \overline{\ll}_3 = 12 e_1 + e_2 + \overline{\ll}_1.
\]
The matrix of relations is 
\[
R_0 = 
\begin{pmatrix}
3 & 0 & -2 & 0 & 0
\\
5 & 0 & 1 & -2 & 0
\\
12 & 1 & 1& 0 & -2
\end{pmatrix}
.
\]
By Definition \ref{def:eq-xi}  and Proposition \ref{prop:mon-inc} we get
\[
\mathcal{I}_3 = (X^{e_1+ e_2}, X^{e_2 +\overline{\ll}_1}, X^{e_1 + \overline{\ll}_3} ) \subset \mathcal{J}_3 =
(X^{e_1+ e_2}, X^{e_2 +{\ll}_1}, X^{e_1 + {\ll}_3} ) \subset \mathcal{O}_Z.
\]
By Corollary \ref{th:inclusion1}, if $p$ is a point in the torus of $Z^\Gamma$,
 then $\mathrm{Jac}(S_p) \mathcal{O}_Z$ is the monomial ideal $X^{\g_0} \mathcal{J}_3$. 
By Remark \ref{rem:pzero},  if $p= 0$ then $\mathrm{Jac}(S_0) \mathcal{O}_Z$ is the monomial ideal $X^{\g_0} \mathcal{I}_3$. 
\end{Exam}

\section{An application} \label{sec:app}

In this section we give a method to compute the number $\overline{\nu}_{\mathrm{Jac}(S)} ( \phi )$ 
associated to a function $\phi  \in  \mathcal{O}_S$ and the Jacobian ideal $\mathrm{Jac}(S)$ of $S$, when $S$ is
an irreducible quasi-ordinary hypersurface of dimension two.

We recall first the definition of the numbers
 $\overline{\nu}_I (a)$ and 
$\overline{\nu}_\mathcal{I} (\mathcal{J})$ associated to an element $a \in A$ and a pair of proper ideals $\mathcal{I}, \mathcal{J}$ of a ring $A$ by Lejeune-Jalabert and Teissier (see \cite{LT}).
These numbers have many interesting applications in the frame of
analytic geometry, in particular in connections with the
 \L
ojasiewicz exponents and analytic arcs (see \cite{LT,
Teissier-inv, Hickel}). 
We show then how to describe the invariant  $\overline{\nu}_\mathcal{I} (\mathcal{J})$ 
explicitly when $A = \C\{ \Lambda \}$ is the ring of germs of functions 
at the origin of an
affine toric variety $Z^\Lambda$ and $\mathcal{I}$ and $\mathcal{J}$ are a pair of
monomial ideals of $A$.

\subsection{The number $\overline{\nu}_I (a)$}
\label{tor-inv}

Let $A$ be a ring and $\mathcal{I} \ne A$ be an ideal. Let us consider first the order function 
$\nu_\mathcal{I}: A \to \Z$
associated with the
$\mathcal{I}$-adic filtration $A =\mathcal{I}^0 \supset \mathcal{I} \supset \mathcal{I}^2 \supset
\cdots$ of $A$, which is defined by 
\[
\nu_\mathcal{I} (a) := \sup \{ n \mid a \in \mathcal{I}^n \} \in \Z_{\geq 0}.\]
If $\mathcal{J}$ is an ideal of $A$ we set $\nu_\mathcal{I} (\mathcal{J})
:= \sup \{ n \mid \mathcal{J} \subset \mathcal{I}^n \}$. 
We have also an order function $\overline{\nu}_\mathcal{I} : A \to \R_{\geq 0} \cup \{ + \infty \}$ defined by 
\[
\overline{\nu}_I (a) := \lim_{k\to \infty} \frac{\nu_I (a^k)}{k} .
\] 
If $\mathcal{J}$ is an ideal of $A$ then we set 
 $\overline{\nu}_\mathcal{I}
(\mathcal{J}) := \lim_{k\to \infty} \frac{\nu_I (\mathcal{J}^k) }{k}$. 
 If $\mathcal{J}$ is generated by $a_1, \dots, a_s$ then we have $\overline{\nu}_\mathcal{I} (
 \mathcal{J}) := \min_{i=1,
\dots, s} \overline{\nu}_I (a_i)$. These numbers belong to $\R_{\geq 0}
\cup \{ \infty \}$. See \cite{LT}.


  \subsection{Normalized blow up of a monomial ideal in a toric variety}

The {\it Newton
polyhedron} of $\phi = \sum c_\a X^\a  \in \C \{ \s^\vee \cap M \}$ is 
the convex hull   ${\mathcal N} (\phi)$
of the union of the  Minkowski sum  of sets ${\a} + \s^\vee$, for $\a$ such that $c_\a\ne 0$. 
Similarly, the {\it Newton  polyhedron} of a monomial ideal $\mathcal{I} = ({X}^{u} / u \in {I})$ of $\C\{ \s^\vee \cap M \}$  is the convex hull   ${\mathcal N} (\mathcal{I})$
of the Minkowski sum  of sets ${I} + \s^\vee$.
If $\mathcal{N}$ is the Newton polyhedron of a function or a monomial ideal we
denote by $\mbox{\rm ord}_{\mathcal{N}}$ the support
function of the polyhedron ${\mathcal N}$,
which is defined by
 \[\mbox{\rm ord}_{\mathcal{N} }: \s \rightarrow
\R, \quad  \nu \mapsto \inf_{\omega \in {\mathcal N} } \langle \nu, \omega \rangle.\] 

 The {\it dual fan}   $\Sigma ({\mathcal
I})$ associated to the  polyhedron ${\mathcal N} (\mathcal{I})$ 
consists of  the cones
\[
\s( {\mathcal F} ) := \{  \eta \in \s \; \mid \langle \eta , v
\rangle   = \mbox{\rm ord}_{I} ( \eta ),
 \; \forall v \in {\mathcal F}\},
 \]
for ${\mathcal F}$ running through the faces of ${\mathcal N}
(\mathcal{I})$.


 If $\Sigma = \Sigma (\mathcal{I}) $, the associated the  toric modification
 $
 \p_\Sigma : Z_\Sigma \rightarrow Z_\s$
 is the {\it normalized blow up} of $Z_\s$ centered at the
 monomial ideal $\mathcal{I}$ of  $\C \{ \s^\vee \cap M \}$ (see
 \cite[Chapter 1]{TE} or \cite{LR}).

Let $\{ D_{k} \}_{k \in K(\mathcal{I})}$ be the irreducible components
of the exceptional divisor of $\pi_\Sigma$ which is defined by
 the ideal sheaf  $I \mathcal{O}_{Z_{\Sigma}}$. Each
component $D_k$ is an invariant divisor for the torus action and the
divisorial valuation $\nu_{D_k}$  defined by $D_k$  is a monomial
valuation of the form 
\[ \nu_{D_k} (\sum  a_m X^m ) := \inf_{a_m \ne
0} \langle n_k, m \rangle,\] 
where  $n_k$ is a primitive vector for
the lattice $N$ in an edge $\r$ of $\Sigma$ such that, if  $\t$ is
the unique face of $ \s$ with $\r \subset
\mathrm{int}({\t})$, then the orbit $\O_\t$ is contained in the zero
locus of $\mathcal{I}$ in $Z^\Lambda$. Then, we have
\begin{equation}\label{mult-nb}
{\nu_{D_k} (\mathcal{I}\mathcal{O}_{Z_\Sigma} ) } 
= \min_{u \in I} \nu_{D_k} (X^u) 
= \mbox{\rm ord}_{\mathcal{I}} (n_k).
\end{equation}
See \cite{LR, TE, Grenoble} for more details.

\subsection{The number $\nu_I (\phi)$ in the toric case}

Suppose from now on that $A = \C\{\Lambda\}$ and $\mathcal{I} = ({X}^{u} \mid u \in {I})$, $\phi \in \C\{\Lambda\}$ and
$\mathcal{J}  =  ({X}^{u} \mid u \in {J})$ are
proper monomial ideals of $\C\{\Lambda\}$, for $\Lambda$ as in Section
\ref{sec-tor}. By  \cite[Prop. 0.20]{LT} we have 
\[
\overline{\nu}_\mathcal{I} (\phi) = \overline{\nu}_{\mathcal{I} \C \{ \s^\vee \cap M \} } (\phi)
\mbox{ and } \overline{\nu}_\mathcal{I} (\mathcal{J} ) = \overline{\nu}_{\mathcal{I} \C \{ \s^\vee \cap M \} } ( \mathcal{J}
\C\{ \s^\vee \cap M \} ).
\] 
In order to describe the numbers  $\overline{\nu}_\mathcal{I} (\phi ) $ and $\overline{\nu}_\mathcal{I} (\mathcal{J} ) $ 
we can assume that $\Lambda =
\s^\vee \cap M$.


\begin{Pro} \label{analytic}
With the previous notations we have
\[  
\overline{\nu}_{\mathcal{I}} (\phi) = \min_{k \in K(\mathcal{I})} \frac{ \mbox{\rm ord}_{\mathcal{\mathcal{N} (\phi)}}
(n_k) }{ \mbox{\rm ord}_{\mathcal{I}} (n_k) }
\quad \mbox{ and } \quad
\overline{\nu}_{\mathcal{I}} (\mathcal{J}) = \min_{k \in K(\mathcal{I})} \frac{ \mbox{\rm ord}_{\mathcal{J}}
(n_k) }{ \mbox{\rm ord}_{\mathcal{I}} (n_k) }. \]
\end{Pro}
\begin{proof} If $k \in K(\mathcal{I})$ then we have $ {\nu_{D_k} (\mathcal{I}
\mathcal{O}_{Z_\Sigma} ) } = \mbox{\rm ord}_{\mathcal{I}} (n_k)$ by
\eqref{mult-nb}.
By the results of \S 4.1 in \cite{LT} 
if $\f \in \C\{\s^\vee \cap M \}$ then the following equality holds
\begin{equation} \label{eq:appendix1}
\overline{\nu}_\mathcal{I} (\f)  = \min_{k
\in K(\mathcal{I})} \frac{\nu_{D_k} (\f)}{\nu_{D_k} ( \mathcal{I} \mathcal{O}_{Z_\Sigma}
) }.
\end{equation}
Since $\nu_{D_k}$ is the monomial valuation defined by the primitive integral vector $n_k$, we get
$
\nu_{D_k} (\f) = \ord_{\mathcal{N}(\phi)} (n_k)
$.

If $\mathcal{J} =
(X^{m_1}, \dots, X^{m_r})$ we have $\nu_{D_k} (X^{m_i}) =
\langle n_k, m_i \rangle$.   We get
\[  
\begin{array}{lcl}
\overline{\nu}_{\mathcal{I}}
(\mathcal{J}) & =  & \min_{i =1, \dots,r} \min_{k \in K(\mathcal{I})} \langle n_k, m_i
\rangle  ( \mbox{\rm ord}_{\mathcal{I}} (n_k) )^{-1} 
\\ 
& = &
\min_{k \in K(\mathcal{I})} \min_{i =1, \dots, r} \langle n_k, m_i \rangle  ( \mbox{\rm ord}_{\mathcal{I}} (n_k) )^{-1}
\\ 
& = &
 \min_{k \in K(\mathcal{I})} \frac{ \mbox{\rm ord}_{\mathcal{J}}
(n_k) }{ \mbox{\rm ord}_{\mathcal{I}} (n_k) }.
\end{array}
\]
\end{proof}

\subsection{Application to the case of a quasi-ordinary hypersurface $S$ of dimension two}

Let $(S,0) \subset (\C^3,0)$ denote an irreducible germ of quasi-ordinary hypersurface of dimension two parametrized by a
normalized quasi-ordinary branch $\zeta$. Denote by $Z $ the normalization of $S$. 
It is the germ of affine toric variety $Z^{\s^\vee \cap M}$ at its origin. We use the notations of Section \ref{qo}. 

By Corollary \ref{cor:inclusion1} and Theorem \ref{jotak2}  the ideal 
$\mathrm{Jac}(S) \mathcal{O}_Z = X^{\g_0} \mathcal{J}_g$ is a monomial ideal 
determined by the characteristic exponents of $\z$. 
By  \cite[Prop. 0.20]{LT}  we have 
$\overline{\nu}_{\mathrm{Jac}(S)} ( \phi )  =  \overline{\nu}_{\mathrm{Jac}(S) \mathcal{O}_Z} ( \phi ) $.
Then, we can apply 
Proposition \ref{analytic} to determine $\overline{\nu}_{\mathrm{Jac}(S)} ( \phi )$ from the 
Newton polygons of $\phi$ and of the monomial ideal 
$\mathcal{J}_g$
of 
$\C\{ \s^\vee \cap M\}$.

\subsection*{Acknowledgements}
The author thanks Evelia Garc\'{\i}a Barroso, Antoni Rangachev and Bernard Teissier and the referee  for their comments, which helped to improve the presentation of the paper.  This work is supported by the Spanish grants of MCIN (PID2024-156181NB-C32 and PID2020-114750GB-C32/AEI/10.13039/501100011033). 

\normalsize
\renewcommand{\MR}[1]{}

 \end{document}